\documentclass[12pt]{amsart}
\usepackage{eufrak}
\usepackage[total={6.5in,8.75in}, top=1.2in, left=1.0in, includefoot]{geometry}
\usepackage{ amssymb }
\usepackage{graphics}
\usepackage{ mathrsfs }

\theoremstyle{definition}

\theoremstyle{remark}

\numberwithin{equation}{section}

\begin{document}

\title{Monodromy of Galois representations and equal-rank subalgebra equivalence}

\author{Chun Yin Hui}
 \address{Einstein Institute of Mathematics \\
Edmond J. Safra Campus, Givat Ram \\
The Hebrew University of Jerusalem \\
Jerusalem, 91904, Israel}

  \email{chhui@indiana.edu, pslnfq@gmail.com}

\begin{abstract}
Let $K$ be a number field, $\mathscr{P}$ the set of prime numbers, and $\{\rho_\ell\}_{\ell\in \mathscr{P}}$ a compatible system (in the sense of Serre \cite{Serre3}) of semisimple, $n$-dimensional $\ell$-adic representations of $\mathrm{Gal}(\overline{K}/K)$. Denote the Zariski closure of $\rho_\ell(\mathrm{Gal}(\overline{K}/K))$ in $\mathrm{GL}_{n,\mathbb{Q}_\ell}$ by $G_\ell$ and its Lie algebra by $\mathfrak{g}_\ell$. It is known that the identity component $G_\ell^\circ$ is reductive and the formal character of the tautological representation $G_\ell^\circ\hookrightarrow \mathrm{GL}_{n,\mathbb{Q}_\ell}$ is independent of $\ell$ (Serre). We use the theory of abelian $\ell$-adic representations
to prove that the formal character of the tautological representation of the derived group $(G_\ell^\circ)^\mathrm{der}\hookrightarrow \mathrm{GL}_{n,\mathbb{Q}_\ell}$  is likewise independent of $\ell$. By investigating the geometry of weights of this faithful representation, we prove that the semisimple parts of $\mathfrak{g}_\ell\otimes\mathbb{C}$ satisfy an equal-rank subalgebra equivalence for all $\ell$ which is equivalent to the number of $A_n:=\mathfrak{sl}_{n+1,\mathbb{C}}$ factors for $n\in\{6,9,10,11,...\}$ and the parity of the number of $A_4$ factors in $\mathfrak{g}_\ell\otimes\mathbb{C}$ are independent of $\ell$.

\end{abstract}

\maketitle

\begin{center}
\textbf{$\mathsection1$. Introduction}
\end{center}

Let $K$ be a number field, $\overline{K}$ its algebraic closure and $G_K:=\mathrm{Gal}(\overline{K}/K)$ the absolute Galois group of $K$. Let $\mathscr{P}$ be the set of prime numbers. A compatible system (Definition 3.4) of $\ell$-adic representations $\{\rho_\ell\}_{\ell\in \mathscr{P}}$ of $G_K$ is a collection of continuous representations 
\begin{equation*}
\rho_\ell: G_K \longrightarrow \mathrm{GL}_n(\mathbb{Q_\ell}),
\end{equation*}
indexed by $\mathscr{P}$, such that for any $\ell\in\mathscr{P}$, $\rho_\ell$ is unramified at all but finitely many non-Archimedean places of $K$, and for any $\ell,\ell'\in\mathscr{P}$, the characteristic polynomials of $\rho_\ell(\mathrm{Frob}_w)$ and $\rho_{\ell'}(\mathrm{Frob}_w)$ (well-defined if $\rho_\ell$ and $\rho_{\ell'}$ are unramified at $v$, and $w$ is a valuation of $\overline{K}$ extending $v$) are equal with rational coefficients for all but finitely many non-Archimedean places $v$ of $K$. Such a compatible system arises, for example, from the Galois action of $G_K$ on the $\ell$-adic Tate modules $T_\ell(A)\otimes\mathbb{Q_\ell}$ of an abelian variety $A$ defined over $K$ (see \cite[Chap. 1]{Serre3}) or more generally, on the $\ell$-adic \'etale cohomology groups $H^k_{et}(X_{\overline{K}},\mathbb{Q_\ell})$ of a complete non-singular variety $X$ defined over $K$ (Deligne \cite{Del}). We may assume our representations $\rho_\ell$ are semisimple for all $\ell$ by semi-simplification since the characteristic polynomials of matrices only carry information from the semisimple parts.\\

For a compatible system $\{\rho_\ell\}_{\ell\in \mathscr{P}}$ of semisimple $\ell$-adic representations of $G_K$,
the Zariski closure of $\rho_\ell(G_K)$ in the algebraic group $\mathrm{GL}_{n,\mathbb{Q_\ell}}$ is a reductive algebraic group and is called the \emph{algebraic monodromy group} at $\ell$. Denote it by $G_\ell$, its connected component by $G_\ell^\circ$,  and its Lie algebra by $\mathfrak{g}_\ell$. Let $\Phi_\ell: G_\ell\hookrightarrow \mathrm{GL}_{n,\mathbb{Q}_\ell}$ be the tautological representation. Consider the following conjectures:\\

\noindent\textbf{Conjecture 1.1.} There exists a faithful representation $\Phi:G\hookrightarrow \mathrm{GL}_{n,\mathbb{Q}}$ of  a reductive group $G$ over $\mathbb{Q}$, such that for all $\ell\in \mathscr{P}$,
$(G_\ell,\Phi_\ell)\cong (G,\Phi)\times_\mathbb{Q}\mathbb{Q}_\ell$.\\

\noindent\textbf{Conjecture 1.2.} There exists a faithful representation $\Phi:G\hookrightarrow \mathrm{GL}_{n,\mathbb{Q}}$ of  a connected reductive group $G$ over $\mathbb{Q}$, such that for all $\ell\in \mathscr{P}$,
$(G_\ell^\circ,\Phi_\ell)\cong (G,\Phi)\times_\mathbb{Q}\mathbb{Q}_\ell$.\\

\noindent Conjecture 1.1 is true if $\rho_\ell(G_K)$ is abelian for one $\ell$, see chapter $3$ of \cite{Serre3}. Weaker variants of this conjecture have appeared many times in the literature. For example, if our compatible system $\{\rho_\ell\}_{\ell\in \mathscr{P}}$ comes from the Galois action on the Tate modules of an abelian variety, then Conjecture 1.2 follows immediately from the semisimplicity of $\rho_\ell$ (Faltings \cite{F}) and the Mumford-Tate conjecture (see Mumford \cite{Mu1}, Serre \cite[$\mathsection9$]{Serre4}). If $\{\rho_\ell\}_{\ell\in \mathscr{P}}$ comes from the Galois action on the \'etale cohomology groups of a complete non-singular variety, then Conjecture 1.2 is implied (see \cite[$\mathsection5$]{LP3}) by the well-known semisimplicity conjecture \cite[$\mathsection9$]{Serre4} and the general Tate conjecture \cite{Tate}. By the method of Frobenius tori, Serre \cite[p.~6,17]{Serre1},\cite[$\mathsection$2.2.3]{Serre2} has proved the following $\ell$-independence results.\\

\noindent\textbf{Theorem 1.3.} The open subgroup of finite index $\rho_\ell^{-1}(G_\ell^\circ(\mathbb{Q_\ell}))\subset G_K$ is independent of $\ell$.\\

\noindent\textbf{Theorem 1.4.} The pair $(T_\ell,\Phi_\ell)$ consisting of a maximal torus $T_\ell$ of $G_\ell^\circ$ and the tautological representation $\Phi_\ell:T_\ell\hookrightarrow\mathrm{GL}_{n,\mathbb{Q}_\ell}$  is independent of $\ell$. Therefore, the formal character of $G_\ell^\circ \hookrightarrow\mathrm{GL}_{n,\mathbb{Q}_\ell}$ and hence the rank of $G_\ell^\circ$ are independent of $\ell$.\\

\noindent\textbf{Remark 1.5.} Serre originally stated these results for representations associated to abelian varieties, but his proofs work for arbitrary compatible systems.\\

\noindent\textbf{Remark 1.6.} Theorem 1.3 implies there is a smallest finite extension $K^\mathrm{conn}$ of $K$ such that the Zariski closure of $\rho_\ell(G_{K^\mathrm{conn}})$ in $\mathrm{GL}_{n,\mathbb{Q_\ell}}$ is equal to $G_\ell^\circ$ for all $\ell$. For more results about $K^\mathrm{conn}$, see Silverberg and Zarhin \cite{SZ1},\cite{SZ2} and Larsen and Pink \cite{LP2}.\\

We study $\ell$-independence of compatible systems of $\ell$-adic representations of $G_K$, assuming semisimplicity. We are not able to prove Conjecture 1.1 and Conjecture 1.2, but using the theory of abelian $\ell$-adic representations, we extend Theorem 1.4 by proving that the formal character of the representation of the derived group $(G_\ell^\circ)^\mathrm{der}\hookrightarrow \mathrm{GL}_{n,\mathbb{Q}_\ell}$ is independent of $\ell$. Using these data, we prove that the semisimple parts of the complexified Lie algebras $\mathfrak{g}_\ell\otimes\mathbb{C}$ satisfy an equal-rank subalgebra equivalence (Definition 2.18). The results of this paper are summarized as follows.
 
\begin{enumerate}
	\item (Proposition 3.18) The dimension of the center of $G_\ell^\circ$ is bounded by $d_{K^{\mathrm{conn}}}$, the common dimension of the Serre groups $S_\mathfrak{m}$ associated to the number field $K^\mathrm{conn}$.\\

	\item (Theorem 3.19) The triple $(((G_\ell^\circ)^{\mathrm{der}}\cap T_\ell)^\circ, T_\ell, \Phi_\ell)$ is independent of $\ell$, where $T_\ell$ is a maximal torus of $G_\ell^\circ$ and $\Phi_\ell$ is the embedding of $T_\ell$ into $\mathrm{GL}_{n,\mathbb{Q}_\ell}$. Therefore, the formal character of the tautological representation $(G_\ell^\circ)^\mathrm{der}\hookrightarrow \mathrm{GL}_{n,\mathbb{Q}_\ell}$ and hence the semisimple rank of $G_\ell^\circ$ are independent of $\ell$.\\
	
	\item (Theorem 3.21) Consider the free abelian group of virtual complex simple Lie algebras which contains semisimple Lie algebras naturally. We divide by the subgroup generated by all expressions $\mathfrak{g}-\mathfrak{h}$ where $\mathfrak{h}\subset\mathfrak{g}$ are equal rank semisimple Lie algebras. The semisimple parts of the complexified Lie algebras $\mathfrak{g}_\ell\otimes\mathbb{C}$ have the same image in this quotient group (satisfy equal-rank subalgebra equivalence in Definition 2.18) for all $\ell$. This is equivalent to the number of $A_n:=\mathfrak{sl}_{n+1,\mathbb{C}}$ factors for $n\in\{6,9,10,11,...\}$ and the parity of the number of $A_4$ factors in $\mathfrak{g}_\ell\otimes\mathbb{C}$ are independent of $\ell$. \\
	
	\item (Theorem 4.1) Let $K$ be a field, finitely generated over $\mathbb{Q}$ and $G_K$ its absolute Galois group. If the system $\{\rho_\ell\}_{\ell\in \mathscr{P}}$ arises from the Galois action of $G_K$ on the Tate modules of an abelian variety $X$ defined over a field $K$, then $(3)$ holds for $\mathfrak{g}_\ell\otimes\mathbb{C}$.
\end{enumerate}

\noindent\textbf{Remark 1.7.} Larsen and Pink \cite{LP92} studied compatible systems of representations of profinite groups that are endowed with a dense subset of ``Frobenius'' elements. They gave an example in \cite[$\mathsection10$]{LP92} that the semisimple rank of $G_\ell$ (the algebraic monodromy group at $\ell$) depends on $\ell$, contrary to $(2)$.  \\

\noindent\textbf{Remark 1.8.} (1), (2), (3) above also hold for any semisimple compatible system of $\lambda$-adic \cite[Chap. 1 $\mathsection2.3$]{Serre3} representations. This will be explained at the end of section 3.\\ 

\noindent\textbf{Remark 1.9.} We study $\ell$-independence of mod $\ell$ Galois representations that arise from \'etale cohomology in a subsequent paper \cite{mod l} and obtain ``mod $\ell$" versions of (2) and (3) when $\ell\gg 0$.\\

\noindent Let us sketch the proofs of our results. We may assume $K=K^\mathrm{conn}$ so that $G_\ell=G_\ell^\circ$. Since the Lie algebra of $[G_\ell,G_\ell]$ is the semisimple part of $\mathfrak{g}_\ell$. Hence, we could study the dimension of the center of $\mathfrak{g}_\ell$ by considering the dimension of the image of following semisimple $\ell$-adic representation
\begin{equation*}
G_K\stackrel{\rho_\ell}{\longrightarrow} G_\ell(\mathbb Q_\ell)\rightarrow G_\ell/[G_\ell,G_\ell](\mathbb Q_\ell)\hookrightarrow \mathrm{GL}_m(\mathbb{Q}_\ell).
\end{equation*}

\noindent This map factors through the quotient $G_K\rightarrow G_K^\mathrm{ab}$, so we get $\Psi_\ell: G_K^\mathrm{ab}\rightarrow \mathrm{GL}_m(\mathbb{Q}_\ell)$. Let $\sum_K$ be the set of finite places of $K$. Since $\Psi_\ell$ is unramified at all except finitely many places $S$ of $K$, if $F_w$ is a Frobenius element of a valuation $w$ of $\overline{K}$ extending a place $v\in\sum_K\setminus S$, then the eigenvalues of $\Psi_\ell(F_w)$ are algebraic numbers. This is the key observation which leads to result $(1)$. By imitating the proof of Theorem 1.10, we can find an integer $N$ such that $\Psi_\ell^N:G_K^\mathrm{ab}\rightarrow \mathrm{GL}_m(\mathbb{Q}_\ell)$ is locally algebraic. \\

\noindent\textbf{Theorem 1.10.} (See Serre \cite[Chap. 3 $\mathsection$3]{Serre3}, Waldschmidt \cite{W}, Henniart \cite[$\mathsection$5]{Hen} )\\
\noindent If $\rho: G_K^\mathrm{ab}\rightarrow \mathrm{GL}_m(\mathbb{Q}_\ell)$ is a rational, semisimple, $\ell$-adic abelian representation of $K$, then $\rho$ is locally algebraic.\\
 
\noindent Therefore, by the theory of abelian $\ell$-adic representation, $\Psi_\ell^N$ arises from some abelian $\ell$-adic representation attached to some Serre group $S_\mathfrak{m}$. In other words, $S_\mathfrak{m}(\mathbb{Q}_\ell)$ surjects onto $\Psi_\ell^N(G_K^\mathrm{ab})$. Since the dimension of $\Psi_\ell^N(G_K^\mathrm{ab})$, $\Psi_\ell(G_K^\mathrm{ab})$ and the center of $\mathfrak{g}_\ell$ are the same, we get $(1)$. By the above techniques, the upper bound of the dimension of center in $(1)$ and Theorem 1.4, we construct an auxiliary compatible system of representations of $G_K$ to prove $(2)$. The restriction of the formal character to the derived group of $G_\ell^\circ$ are independent of $\ell$ by $(2)$. Therefore, we only need to study the semisimple part of $\mathfrak{g}_\ell$. We have for each $\ell$ a faithful representation of $(\mathfrak{g}_\ell)_{ss}\otimes\mathbb{C}$ which gives the same formal character. So we need to answer the question below.\\

\noindent\textbf{Q:} To what extent is a complex semisimple Lie algebra $\mathfrak{g}$ determined if the formal character of a faithful representation of $\mathfrak{g}$ is given?\\

\noindent Larsen and Pink have answered this question in the case that the representation is irreducible \cite[$\mathsection$4]{LP1}. The difficulty of the question can be illustrated by the following example. We know that $E_7\times A_1$ and $A_4\times A_4$ are subalgebras of maximal rank in $E_8$, so if we restrict a representation of $E_8$ to $E_7\times A_1$ and $A_4\times A_4$, the formal characters of these two representations are the same. We address this question by investigating the geometry of the roots and weights in the formal characters. Actually, we will first prove that the number of $A_n$ factors for $n\in\{6,9,10,11,...\}$ and the parity of the number of $A_4$ factors in $\mathfrak{g}$ are invariants (Theorem 2.14, 2.17). Then, it follows easily that the image of $\mathfrak{g}$ in the quotient group in $(3)$ is invariant (Theorem 2.19). This together with Proposition 2.20 imply $(3)$. Finally, $(4)$ is a direct consequence of $(3)$ and a result on $\ell$-independence of specialization of monodromy groups of abelian varieties (Hui \cite[Thm. 2.5]{Hui}).\\

The structure of this paper is as follows. Section 2 is devoted to answering question \textbf{Q}, which is purely representation theoretic; we only use the representation theory of complex semisimple Lie algebras. We will prove $(1),(2)$, and $(3)$ in section 3. This section relies heavily on the theory of abelian $\ell$-adic representation (Serre \cite{Serre3}). In section 4, we consider systems $\{\rho_\ell\}_{\ell\in \mathscr{P}}$ coming from abelian varieties and prove $(4)$.\\

\begin{center}
\textbf{$\mathsection2$ Geometry of weights in formal characters}
\end{center}

\noindent This section is purely representation theoretic and self-contained. It is devoted to answering the question below.\\

\noindent\textbf{Q:} To what extent is a complex semisimple Lie algebra $\mathfrak{g}$ determined if the formal character of a faithful representation of $\mathfrak{g}$ is given?\\

\noindent(\textbf{2.1}) Let $\mathfrak{g}$ be a complex semisimple Lie algebra and $t_\mathfrak{g}$ some Cartan subalgebra of $\mathfrak{g}$. Denote the roots, the weight lattice and the Weyl group by $\Phi_\mathfrak{g}$, $\Lambda_\mathfrak{g}$ and $W_\mathfrak{g}$ respectively. If $\Theta:\mathfrak{g}\rightarrow \mathfrak{gl}(V)$ is a representation of $\mathfrak{g}$ on some $n$-dimensional complex vector space $V$,
then the action of $t_\mathfrak{g}$ on $V$ can be diagonalized and then we have $n$ weight vectors $\alpha_1,...,\alpha_n\in t_\mathfrak{g}^*$. Let $\mathbb{Z}[\Lambda_\mathfrak{g}]$ be the group ring over $\mathbb{Z}$ generated by the free abelian group $\Lambda_\mathfrak{g}$. We define the \emph{formal character} of $\Theta$ to be $\mathrm{Char}_\Theta(V):=\alpha_1+\cdots+\alpha_n\in \mathbb{Z}[\Lambda_\mathfrak{g}]$. We know that $\mathrm{Char}_\Theta(V)$ is invariant under the Weyl group $W_\mathfrak{g}$. Let $\Psi:\mathfrak{h}\rightarrow \mathfrak{gl}(V')$ be a representation of another complex semisimple $\mathfrak{h}$ on an $n$-dimensional complex vector space $V'$. We say $\mathfrak{h}$ and $\mathfrak{g}$ \emph{have the same formal character} if there is an isomorphism $F$ between $t_\mathfrak{h}^*$ and $t_\mathfrak{g}^*$ such that $F(\mathrm{Char}_\Psi(V))$(defined in an obvious way) is equal to $\mathrm{Char}_\Theta(V)$. We can now state our main theorems of this section.\\

\noindent \textbf{Theorem 2.14} If faithful representations of two complex semisimple Lie algebras $\mathfrak{g}$ and $\mathfrak{h}$ have the same formal character, then the number of $A_n$ factors of $\mathfrak{g}$ and $\mathfrak{h}$ are the same when $n\in\{6,9,10,...\}$.\\

\noindent \textbf{Theorem 2.17} If faithful representations of two complex semisimple Lie algebras $\mathfrak{g}$ and $\mathfrak{h}$ have the same formal character, then the parities of the numbers of $A_4$ factors of $\mathfrak{g}$ and $\mathfrak{h}$ are the same.\\

\noindent \textbf{Theorem 2.19} If faithful representations of two complex semisimple Lie algebras $\mathfrak{g}$ and $\mathfrak{h}$ have the same formal character, then $\mathfrak{g}$ and $\mathfrak{h}$ satisfy equal-rank subalgebra equivalence (Definition 2.18).\\

\noindent(\textbf{2.2}) If $\Theta:\mathfrak{g}\rightarrow \mathfrak{gl}(V)$ is a faithful representation, then $\mathfrak{g}$ is embedded as a subspace in End$(V)=V\otimes V^*$ through $\Theta$. It is easy to check the representation of $\mathfrak{g}$ on this subspace $\mathfrak{g}\subset \mathrm{End}(V)=V\otimes V^*$ is the adjoint representation \cite[Chap. 13.1]{FH}. Therefore, the representation of $\mathfrak{g}$ on $V\otimes V^*$ contains the adjoint representation of $\mathfrak{g}$ as a subrepresentation and is also faithful. The formal character of $V\otimes V^*$ depends only on the formal character of $V$ and is just the sum of all the differences of weights in $\mathrm{Char}_\Theta(V)$. All the roots of $\mathfrak{g}$ appear in the formal character $\mathrm{Char}(V\otimes V^*)$ because the adjoint representation is a subrepresentation. So from now on, we may further assume that $\mathrm{Char}_\Theta(V)$ contains all the roots of $\mathfrak{g}$. The advantage of this assumption is that there are  strong geometric connections among roots and weights once we introduce a suitable Euclidean metric on $\Lambda_\mathfrak{g}\otimes\mathbb{R}$.\\

\noindent(\textbf{2.3}) Suppose $\mathrm{Char}_\Theta(V)=\alpha_1+\cdots+\alpha_n$, we define an inner product on $(\Lambda_\mathfrak{g}\otimes\mathbb{R})^*$ in terms of the formal character by setting 
\begin{equation*}
(x_1,x_2)=\sum_{i=1}^n\alpha_i(x_1)\alpha_i(x_2).
\end{equation*}
We denote the dual inner product on $\Lambda_\mathfrak{g}\otimes\mathbb{R}$ by $\left\langle~, ~\right\rangle$. Since $\Theta$ is faithful, $\{\alpha_i\}$ spans $\Lambda_\mathfrak{g}\otimes\mathbb{R}$, so $(~,~)$ and $\left\langle~, ~\right\rangle$ are positive definite. Since $\mathrm{Char}_\Theta(V)$ is $W_\mathfrak{g}$ invariant, so is $\left\langle~, ~\right\rangle$. This determines $\left\langle~, ~\right\rangle$ up to a positive scalar factor on each simple root system of $\mathfrak{g}$. $W_\mathfrak{g}$ is then a subgroup of the orthogonal group $\mathrm{O}(\Lambda_\mathfrak{g}\otimes\mathbb{R})$ under this Euclidean inner product. Note that $\left\langle~, ~\right\rangle$ is defined solely by the formal character $\mathrm{Char}_\Theta(V)$.\\

\noindent(\textbf{2.4}) Now if $V$ and $V'$ are faithful representations of $\mathfrak{g}$ and $\mathfrak{h}$ with the same formal character, then we could assume that the roots of $\mathfrak{g}$ and $\mathfrak{h}$ appear in the formal character by (\textbf{2.2}).  We can define an Euclidean inner product (which depends only on the formal character) on $\Lambda_\mathfrak{g}\otimes\mathbb{R}$ and on $\Lambda_\mathfrak{h}\otimes\mathbb{R}$ respectively by (\textbf{2.3}). Since the formal characters are the same, $\Lambda_\mathfrak{g}\otimes\mathbb{R}$ and $\Lambda_\mathfrak{h}\otimes\mathbb{R}$ are isometric.\\

\noindent Let $v\in\Phi_\mathfrak{g}$ and $u\in\Phi_\mathfrak{h}$. We claim that the angle $\theta$ between them with respect to the Euclidean metric above belongs to the set $\{0^\circ, 30^\circ, 45^\circ, 60^\circ, 90^\circ, 120^\circ, 135^\circ, 150^\circ, 180^\circ\}$. Indeed, since $u$ (a root of $\mathfrak{h}$) is also a weight of $\mathfrak{g}$ while $v$ (a root of $\mathfrak{g}$) is also a weight of $\mathfrak{h}$, we still have the following relations \cite[$\mathsection$14.1]{FH}
\begin{equation*}
2\frac{||u||\cos\theta}{||v||}=2\frac{\left\langle u, v\right\rangle}{\left\langle v, v\right\rangle}\in\mathbb{Z},
\end{equation*}
\begin{equation*}
2\frac{||v||\cos\theta}{||u||}=2\frac{\left\langle u, v\right\rangle}{\left\langle u, u\right\rangle}\in\mathbb{Z}.
\end{equation*}
\noindent The product of the left hand side is $4\cos^2\theta\in\mathbb{Z}$ and we obtain our claim.\\

\noindent Now, we can determine the ratio $\frac{||u||}{||v||}$.\\

If $\theta=0^\circ$ or $180^\circ$, then $\frac{||u||}{||v||}$, $\frac{||v||}{||u||}\in\frac{\mathbb{Z}}{2}$. We conclude that   $\frac{||v||}{||u||}\in\{\frac{1}{2}, 1, 2\}$;\\

If $\theta=30^\circ$ or $150^\circ$, then $\frac{||u||}{||v||}$, $\frac{||v||}{||u||}\in\frac{\mathbb{Z}}{\sqrt{3}}$. We conclude that   $\frac{||v||}{||u||}\in\{\frac{1}{\sqrt{3}}, \sqrt{3}\}$;\\

If $\theta=45^\circ$ or $135^\circ$, then $\frac{||u||}{||v||}$, $\frac{||v||}{||u||}\in\frac{\mathbb{Z}}{\sqrt{2}}$. We conclude that    $\frac{||v||}{||u||}\in\{\frac{1}{\sqrt{2}}, \sqrt{2}\}$;\\

If $\theta=60^\circ$ or $120^\circ$, then $\frac{||u||}{||v||}$, $\frac{||v||}{||u||}\in\mathbb{Z}$. We conclude that    $\frac{||v||}{||u||}=1$.\\

\noindent\textbf{Lemma 2.5.} Every semisimple Lie algebra $\mathfrak{g}$ contains a semisimple subalgebra $\mathfrak{g}'$ of maximal rank such that every simple factor of $\mathfrak{g}'$ is of type $A_n$ and the number of $A_n$ factors of $\mathfrak{g}$ and $\mathfrak{g}'$ are the same for $n\in\{4,6,9,10,11,12,...\}$.\\

\noindent\textbf{Proof.} By \textbf{Table $1$} \cite[Table 5]{GOV} at the end of this section and the following facts \cite[Chap. 18, 19]{FH}:
\[ \begin{array}{lcl}

\mathfrak{so}(3)=A_1,\\
\mathfrak{so}(4)=A_1\times A_1,\\
\mathfrak{so}(6)=A_3,\\
\mathfrak{sp}(2)=A_1,\\
\mathfrak{so}(4)\subset \mathfrak{so}(5)=\mathfrak{sp}(4),

\end{array}
\]
 one can choose a subalgebra $\mathfrak{g}'$ of $\mathfrak{g}$ of maximal rank such that every simple factor of $\mathfrak{g}'$ is of type $A_n$ and that the number of $A_n$ factors of $\mathfrak{g}$ and $\mathfrak{g}'$ are the same for $n\in\{4,6,9,10,11,12,...\}$.\qed\\
 
\noindent Since $\mathfrak{g}'\subset\mathfrak{g}$ are of same rank, one may assume they have the same Cartan subalgebra. If $\Theta: \mathfrak{g}\rightarrow \mathfrak{gl}(V)$ is a representation of $\mathfrak{g}$ and we take the restriction to $\mathfrak{g}'$, then $\mathrm{Char}_\Theta(V)=\mathrm{Char}_{\Theta|_{\mathfrak{g}'}}(V)$. To prove Theorem 2.14 and Theorem 2.17, we can reduce to the case that every simple factor of our Lie algebras $\mathfrak{g},\mathfrak{h}$ is of type $A_n$.\\

\noindent\textbf{Lemma 2.6.} Suppose $\mathfrak{q}$ is a simple factor of $\mathfrak{g}$ of type $A_n$, $n\geq 2$. Let $u$ be a root of $\mathfrak{h}$ such that $u\notin (\Lambda_\mathfrak{q}\otimes\mathbb{R})\cup (\Lambda_\mathfrak{q}\otimes\mathbb{R})^\bot\subset \Lambda_\mathfrak{g}\otimes\mathbb{R}$. Then the angle $\theta$ between $u$ and any root of $\mathfrak{g}$ can only be $45^\circ$, $60^\circ$, $90^\circ$, $120^\circ$, $135^\circ$.\\

\noindent\textbf{Proof.} We decompose the semisimple Lie algebra $\mathfrak{g}$ into simple factors with $\mathfrak{q}_1=\mathfrak{q}$.
\begin{equation*}
\Lambda_\mathfrak{g}\otimes\mathbb{R}= (\Lambda_{\mathfrak{q}_1}\otimes\mathbb{R})\bot (\Lambda_{\mathfrak{q}_2}\otimes\mathbb{R})\bot\cdots\bot (\Lambda_{\mathfrak{q}_m}\otimes\mathbb{R})
\end{equation*}
\noindent Since $u\notin (\Lambda_\mathfrak{q}\otimes\mathbb{R})\cup (\Lambda_\mathfrak{q}\otimes\mathbb{R})^\bot$, $u\notin \Lambda_{\mathfrak{q}_i}\otimes\mathbb{R}$ for all $i$. Thus the angle cannot be $0^\circ$ or $180^\circ$. Assume $v$ is a root of $\mathfrak{g}$ such that the angle $\theta$ between $u$ and $v$ is $30^\circ$ (if $\theta=150^\circ$, then choose root $-v$). Since $\mathfrak{q}=\mathfrak{q}_1$ is of type $A_n$ ($n\geq 2$), we can always choose some root $w\in\Phi_{\mathfrak{q}_1}$ such that $\{u, v, w\}$ are linearly independent and the angle between $u$ and $w$ is not $90^\circ$. Now, consider the group $G$ generated by $R_u$, $R_v$, and $R_w$, the three reflections of hyperplanes corresponding to roots $u$, $v$, $w$. Since each reflection is either in the Weyl group of $\mathfrak{g}$ or in the Weyl group of $\mathfrak{h}$, the formal character $\mathrm{Char}_\Theta(V)$ is invariant under $G$. It follows that $G$ is finite because $G$ permutes weights in $\mathrm{Char}_\Theta(V)$ and those weights span $\Lambda_\mathfrak{g}\otimes\mathbb{R}$. So we conclude that $G$ is a finite subgroup of $\mathrm{O}(\Lambda_\mathfrak{g}\otimes\mathbb{R})$, the orthogonal group of $\Lambda_\mathfrak{g}\otimes\mathbb{R}$. Also note that the three-dimensional space spanned by $\{u,v,w\}$ is actually invariant under $G$. Denote the determinant $1$ subgroup of $G$ by $G^+$, it is a finite subgroup of SO$(3)$. By the classification of finite subgroups of SO$(3)$ \cite[Chap. 5 Thm. 9.1]{A}, $G^+$ is one of the following:
\begin{enumerate}
\item[$C_k$]: The cyclic group of rotations by multiples of $360^\circ/k$ about a line;
\item[$D_k$]: The dihedral group of symmetries of regular $k$-gon;
\item[$T$]: The tetrahedral group of $12$ rotations carrying a regular tetrahedron to itself;
\item[$O$]: The octahedral group of $24$ rotations carrying a cube to itself;
\item[$I$]: The icosahedral group of $60$ rotations carrying a regular icosahedron to itself.
\end{enumerate}

Note that both the angle between $u,v$ and the angle between $u,w$ are not $90^\circ$, and this means $G^+$ can only be $T$, $O$ or $I$. But the angle between $u,v$ is $30^\circ$. This implies $G^+$ contains $R_u\circ R_v$, a rotation of order $6$ which is impossible. This finishes the proof.\qed\\

\noindent(\textbf{2.7}) Let $\mathfrak{q}$ be a complex semisimple Lie algebra of type $A_n$, $n\geq1$, then $\Lambda_\mathfrak{q}\otimes_\mathbb{Z}\mathbb{R}$ has dimension $n$. There exist weights $e_1,...,e_n,e_{n+1}\in \Lambda_\mathfrak{q}$ such that (See \cite[Table 1]{GOV})
\begin{equation*}
\Lambda_\mathfrak{q}\otimes\mathbb{R}=\mathrm{span}_\mathbb{R}\{e_1,...,e_n,e_{n+1}\},
\end{equation*}
\begin{equation*}
e_1+e_2+...+e_n+e_{n+1}=0,
\end{equation*}
\begin{equation*}
\Lambda_\mathfrak{q}=\mathbb{Z}e_1\oplus\cdots\oplus\mathbb{Z}e_n.
\end{equation*}
\noindent If we normalize so that the length of roots is $\sqrt{2}$, then 
\begin{equation*}
\left\langle e_i, e_j \right\rangle= \left\{ \begin{array}{lll}
 \frac{n}{n+1} &\mbox{if}& 1\leq i=j\leq n\\
 \frac{-1}{n+1} &\mbox{if}& 1\leq i\neq j\leq n
\end{array}\right.
\end{equation*}

\noindent The roots $\Phi_\mathfrak{q}$ comprise the set $\{e_i-e_j:$ $1\leq i\neq j\leq n+1\}$. We can take $\{e_1,...,e_n\}$ as a basis of $\Lambda_\mathfrak{q}\otimes\mathbb{R}$.The picture for $A_2$ looks like\\

\begin{center}\scalebox{.5}{\includegraphics{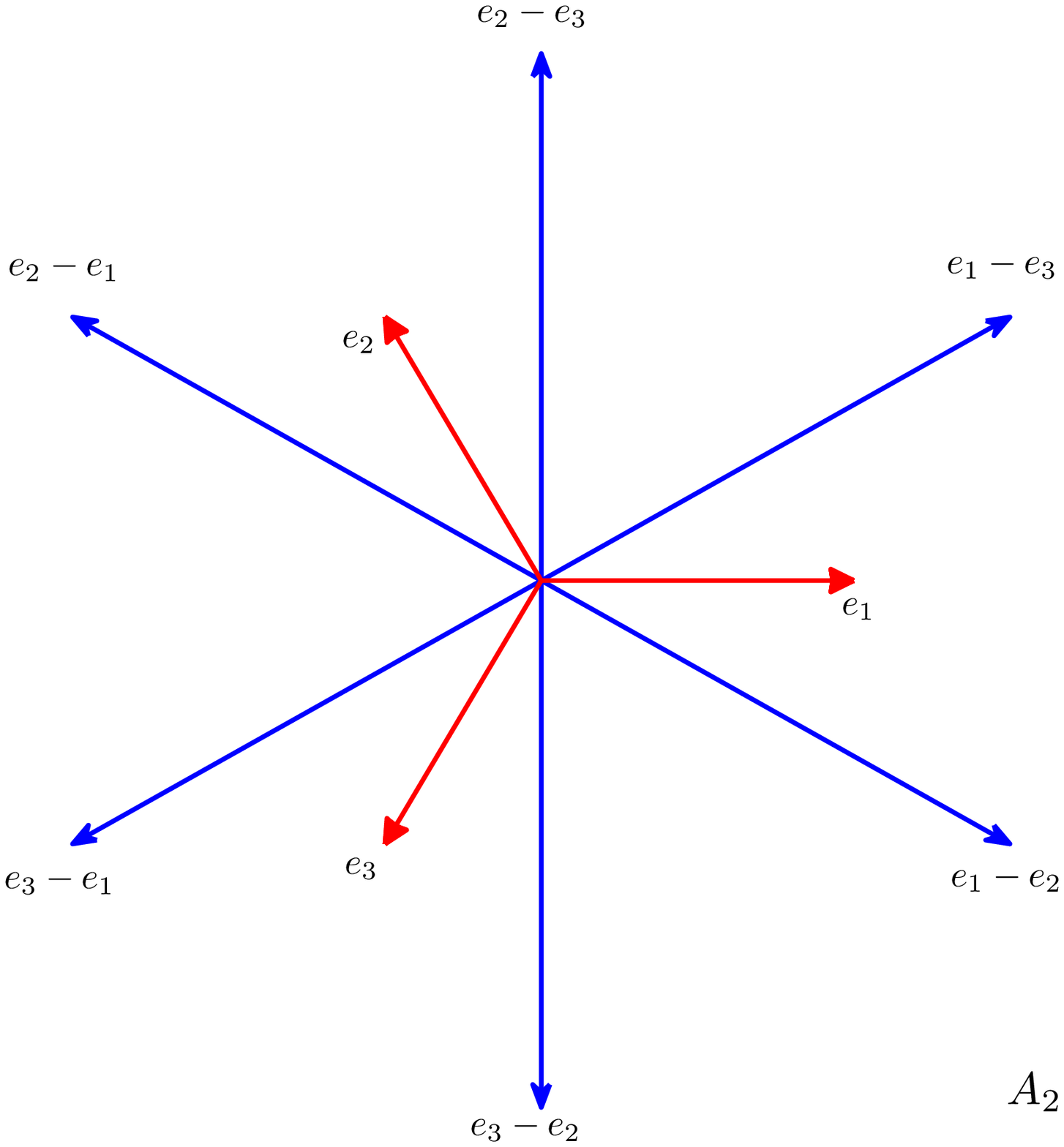}}\end{center}

\noindent The square of the length of a weight $w=a_1e_1+\cdots+a_ne_n$ is
\begin{equation*}
\left\langle a_1e_1+\cdots+a_ne_n, a_1e_1+\cdots+a_ne_n \right\rangle=\frac{\sum_{i=1}^nna_i^2-2\sum_{i<j}a_ia_j}{n+1}
\end{equation*}
\begin{equation}
=\frac{\sum_{i=1}^na_i^2+\sum_{i< j}(a_i-a_j)^2}{n+1}\geq \frac{(n-k)+k(n-k)}{n+1}, \tag{$\Omega$}
\end{equation}
if $k$ is the number of zero $a_i$. Thus if $w\neq 0$ satisfies $||w||<1$, then
\begin{equation*}
1> \frac{(n-k)+k(n-k)}{n+1},
\end{equation*}
\noindent so
\begin{equation*}
1> k(n-k-1),
\end{equation*}
which implies $k=0$ or $n-1$. Thus, one sees easily that the shortest weight is of length $\sqrt{\frac{n}{n+1}}$. Since $n$ (the``numerator") is the biggest integer smaller than $n+1$ (the``denominator"), it is easy to check that $\sqrt{\frac{n}{n+1}}$ is the only length of weight that is less than $1$.\\

\noindent If $||w||=1$, then it is not hard to see by investigating inequality ($\Omega$) that 
\begin{equation*}
||w||= \frac{(n-k)+k(n-k)}{n+1}.
\end{equation*}
\noindent So
\begin{equation*}
1= k(n-k-1),
\end{equation*}
which implies $k=1$ and $n=3$. \\

\noindent\textbf{Proposition 2.8.} Suppose $\mathfrak{q}$ is a simple factor of $\mathfrak{g}$ of type $A_n$, normalized as in (\textbf{2.7}). Let $u$ be a root of $\mathfrak{h}$. The orthogonal projection of $u$ to $\Lambda_\mathfrak{q}\otimes\mathbb{R}$, denoted by $u'$, belongs to $\Lambda_\mathfrak{q}$. We write $u'=a_1e_1+\cdots+a_ne_n$. If $u\notin (\Lambda_\mathfrak{q}\otimes\mathbb{R})\cup (\Lambda_\mathfrak{q}\otimes\mathbb{R})^\bot\subset \Lambda_\mathfrak{g}\otimes\mathbb{R}$, then the following are true.
\begin{enumerate}
	\item If the angle $\theta$ between $u$ and some root of $\mathfrak{q}$ is $60^\circ$, then either all $a_i\in\{0,1\}$ or all $a_i\in\{0,-1\}$.
	\item If the angle $\theta$ between $u$ and some root of $\mathfrak{q}$ is $45^\circ$, then all $a_i\in\{0,1\}$, or all $a_i\in\{0,-1\}$, or all $a_i\in\{0,2\}$, or all $a_i\in\{0,-2\}$.\\
\end{enumerate}

\noindent\textbf{Proof.} Suppose $\mathfrak{g}=\mathfrak{q}_1\oplus\cdots\oplus \mathfrak{q}_m$ is the decomposition of $\mathfrak{g}$ into a direct sum of simple factors with $\mathfrak{q}_1=\mathfrak{q}$. We know that the representation ring of $\mathfrak{g}$ is the tensor product over $\mathbb{Z}$ of the representation rings of $\mathfrak{q}_i$ for all $i$ \cite[Exercise 23.42]{FH}. In other words, any representation of $\mathfrak{g}$ is given by the direct sum of some tensor products of irreducible representation of $\mathfrak{q}_i$. Therefore, if $u$ is a weight of $\mathfrak{g}$ appearing in $\mathrm{Char}_\Theta(V)$, then $u'$ is some weight appearing in some irreducible representation of $\mathfrak{q}$, so $u'$ belongs to $\Lambda_\mathfrak{q}$.\\

\noindent Since $u\notin (\Lambda_\mathfrak{q}\otimes\mathbb{R})\cup (\Lambda_\mathfrak{q}\otimes\mathbb{R})^\bot$ and the roots $\Phi_\mathfrak{q}$ span $\Lambda_\mathfrak{q}\otimes\mathbb{R}$, there exist some root of $\mathfrak{q}$ such that $\theta\in\{30^\circ, 45^\circ, 60^\circ, 120^\circ, 135^\circ, 150^\circ\}$ by (\textbf{2.4}). I claim that if the angle $\theta$ between $u$ and some root of $\mathfrak{q}$ is $60^\circ$, then the angle between $u$ and any root of $\mathfrak{q}$ belongs to $\{60^\circ, 90^\circ, 120^\circ\}$, and if $\theta=45^\circ$, then the angle between $u$ and any root of $\mathfrak{q}$ belongs to $\{45^\circ, 90^\circ, 135^\circ\}$. Suppose not, let $v_1,v_2\in\Phi_\mathfrak{q}$ such that the angle between $u$ and $v_1$ is $45^\circ$ while the angle between $u$ and $v_2$ is $60^\circ$ (WLOG). Since $\mathfrak{q}$ is of type $A_n$, we have
\begin{equation*}
\frac{||u||}{||v_1||}=\frac{||u||}{||v_2||}.
\end{equation*}
\noindent By (\textbf{2.4}), the right hand number and the left hand number belong to $\mathbb{Z}-\{0\}$ and  $\frac{\mathbb{Z}}{\sqrt{2}}-\{0\}$ respectively but the two sets are disjoint. So this is impossible.\\

\noindent Let $v$ be a root of $\mathfrak{q}$ and $\theta$ the angle between $u$ and $v$, and write $u'=a_1e_1+\cdots+a_ne_n$. First consider $v=e_i-e_j$ where $i,j\leq n$. By (\textbf{2.7}), we have 
\begin{equation*}
\left\langle u', v \right\rangle=\left\langle a_1e_1+\cdots +a_ne_n, e_i-e_j \right\rangle=\frac{a_i(n+1)-\sum_{k=1}^na_k}{n+1}-\frac{a_j(n+1)-\sum_{k=1}^na_k}{n+1}=a_i-a_j.
\end{equation*}

\noindent Then consider $v=e_i-e_{n+1}$ where $i\leq n$. By (\textbf{2.7}), we have $v=e_1+e_2+\cdots +e_{i-1}+2e_i+e_{i+1}+\cdots +e_n$, so
\begin{equation*}
\left\langle u', v \right\rangle=\left\langle a_1e_1+\cdots+a_ne_n, e_1+e_2+\cdots+e_{i-1}+2e_i+e_{i+1}+\cdots+e_n \right\rangle
\end{equation*}
\begin{equation*}
=\frac{a_i(n+1)-\sum_{k=1}^na_k}{n+1}+\sum_{j=1}^n\frac{a_j(n+1)-\sum_{k=1}^na_k}{n+1}=a_i.
\end{equation*}

\noindent If we are in the $60^\circ$ case, then $||u||=||v||=\sqrt{2}$. Therefore, $|a_i|$ and $|a_i-a_j|$ are both of the form
\begin{equation*}
|\left\langle u', v \right\rangle|=|2\frac{\left\langle u', v \right\rangle}{\left\langle v, v \right\rangle}|=|2\frac{\left\langle u, v \right\rangle}{\left\langle v, v \right\rangle}|=|2\cos\theta|
= \left\{ \begin{array}{lll}
 1 &\mbox{if}& \theta=60^\circ, 120^\circ\\
 0 &\mbox{if}& \theta=90^\circ
\end{array}\right.
\end{equation*}

\noindent From this, it is easy to see that either all $a_i\in\{0,1\}$ or all $a_i\in\{0,-1\}$.\\

\noindent If we are in the $45^\circ$ case, then $||u||$ is either $2$ or $1$. If $||u||=2$, then

\begin{equation*}
|a_i|,|a_i-a_j|=|\left\langle u', v \right\rangle|=|2\frac{\left\langle u', v \right\rangle}{\left\langle v, v \right\rangle}|=|2\frac{\left\langle u, v \right\rangle}{\left\langle v, v \right\rangle}|=|2\sqrt{2}\cos\theta|
= \left\{ \begin{array}{lll}
 2 &\mbox{if}& \theta=45^\circ, 135^\circ\\
 0 &\mbox{if}& \theta=90^\circ.
\end{array}\right.
\end{equation*}

\noindent When $||u||=1$, then

\begin{equation*}
|a_i|,|a_i-a_j|=|\left\langle u', v \right\rangle|=|2\frac{\left\langle u', v \right\rangle}{\left\langle v, v \right\rangle}|=|2\frac{\left\langle u, v \right\rangle}{\left\langle v, v \right\rangle}|=|\sqrt{2}\cos\theta|
= \left\{ \begin{array}{lll}
 1 &\mbox{if}& \theta=45^\circ, 135^\circ\\
 0 &\mbox{if}& \theta=90^\circ.
\end{array}\right.
\end{equation*}

\noindent From this, it is easy to see either all $a_i\in\{0,1\}$ or all $a_i\in\{0,-1\}$ or all $a_i\in\{0,2\}$ or all $a_i\in\{0,-2\}$.\qed\\

\noindent(\textbf{2.9}) Again, suppose $\mathfrak{q}$ is a simple factor of $\mathfrak{g}$ of type $A_n$. Let $u$ be a root of $\mathfrak{h}$ such that $u\notin (\Lambda_\mathfrak{q}\otimes\mathbb{R})\cup (\Lambda_\mathfrak{q}\otimes\mathbb{R})^\bot\subset \Lambda_\mathfrak{g}\otimes\mathbb{R}$.  Let $\mathfrak{q}'$ be a simple factor of $\mathfrak{g}$ of rank $n'$ such that the orthogonal projection $u'$ of $u$ to $\Lambda_{\mathfrak{q}'}\otimes\mathbb{R}$ is non-zero. Write $u'=a_1e_1+\cdots+a_{n'}e_{n'}$ as in (\textbf{2.7}). Let $k'$ be the number of zero coefficients. Normalize so that $||u||^2=2$. We can use Proposition 2.8 (coordinates computation) to compute $||u'||^2$. Consider the following $3$ cases:\\
\begin{enumerate}
\item If the angle between $u$ and some root of $\mathfrak{q'}$ is $60^\circ$, then by (\textbf{2.4}) the length of roots of $\mathfrak{q}'$ is $\sqrt{2}$. By using (\textbf{2.7}) and Proposition 2.8, we have
\begin{equation*}
||u'||^2=||a_1e_1+\cdots+a_{n'}e_{n'}||^2=\frac{\sum_{i=1}^{n'}a_i^2+\sum_{1\leq i<j\leq n'}(a_i-a_j)^2}{n'+1}=\frac{(n'-k')(k'+1)}{n'+1}.
\end{equation*}
~\\
\item If the angle between $u$ and some root of $\mathfrak{q}'$ is $45^\circ$ and assuming the length of roots of $\mathfrak{q'}$ is $2$ (see (\textbf{2.4})), then we multiply by $2$ and use (\textbf{2.7}) and Proposition 2.8 to get
\begin{equation*}
||u'||^2=2||a_1e_1+\cdots+a_{n'}e_{n'}||^2=\frac{2(\sum_{i=1}^{n'}a_i^2+\sum_{1\leq i<j\leq n'}(a_i-a_j)^2)}{n'+1}=\frac{2(n'-k')(k'+1)}{n'+1}.
\end{equation*}
~\\
\item If the angle between $u$ and some root of $\mathfrak{q}'$ is $45^\circ$ and assuming the length of roots of $\mathfrak{q}'$ is $1$, then we multiply $1/2$ and use (\textbf{2.7}) and Proposition 2.8 to get
\begin{equation*}
||u'||^2=\frac{||a_1e_1+\cdots+a_{n'}e_{n'}||^2}{2}=\frac{\sum_{i=1}^{n'}a_i^2+\sum_{1\leq i<j\leq n'}(a_i-a_j)^2}{2(n'+1)}
\end{equation*}
\begin{equation*}
=\frac{4(n'-k')(k'+1)}{2(n'+1)}=\frac{2(n'-k')(k'+1)}{n'+1}.
\end{equation*}
\end{enumerate}
~\\
\noindent Note that $(n'-k')+(k'+1)=n'+1$, so it is easy to see $||u'||^2\geq\frac{1}{2}$ for the $60^\circ$ case and $||u'||^2\geq 1$ for the $45^\circ$ case.\\

\noindent(\textbf{2.10}) The computations in (\textbf{2.9}) imply strong restrictions on the lengths of various projections of the root $u$. Suppose we are in the situation of (\textbf{2.9}) and $\mathfrak{q}=\mathfrak{q}_1$ is of type $A_{n}$, $n\geq 4$. Consider the situation that there exists another simple factor $\mathfrak{q}_2$ of rank $n_2$ of $\mathfrak{g}$ such that the projection $u_2$ of $u$ to $\Lambda_{\mathfrak{q}_2}\otimes\mathbb{R}$ is non-zero.\\

 \noindent Let's first show that the angle between $u$ and any root of $\mathfrak{q}_2$ cannot be $30^\circ$ or $150^\circ$. Suppose not, choose root $v_2$ of $\mathfrak{q}_2$ that makes angle $30^\circ$ with $u$. Choose a root $v_1$ of $\mathfrak{q}_1$ so that the angle $\theta$ between $v_1$ and $u$ is less than $90^\circ$ by $u\notin (\Lambda_\mathfrak{q}\otimes\mathbb{R})\cup (\Lambda_\mathfrak{q}\otimes\mathbb{R})^\bot$. (\textbf{2.4}) implies $\theta$ is $30^\circ$ or $45^\circ$ or $60^\circ$. By considering the Euclidean space spanned by $u,v_1,v_2$, one sees easily that the first two cases are impossible because $v_1$ and $v_2$ are perpendicular. For the case $\theta=60^\circ$, one observes that $u$ lies on the plane spanned by roots $v_1,v_2$. Since $n\geq 2$ and $u\notin (\Lambda_\mathfrak{q}\otimes\mathbb{R})\cup (\Lambda_\mathfrak{q}\otimes\mathbb{R})^\bot$, the angle between any root of $\mathfrak g_1$ and $u$ belongs to $\{60^\circ, 90^\circ, 120^\circ\}$ by (\textbf{2.4}) and one can choose a root $v_1'$ of $\mathfrak{q}_1$ making an angle of $60^\circ$ with $u$ by the second paragraph of the proof of Proposition 2.8. This is also absurd if one consider the Euclidean $3$-space spanned by $v_1,v_1',v_2$. Therefore, we can apply the results in (\textbf{2.9}).\\

\noindent Let $u_1$ be the projection of $u$ to $\Lambda_{\mathfrak{q}_1}\otimes\mathbb{R}$. Then we have 

\begin{equation*}
2=||u||^2\geq ||u_1||^2+||u_2||^2.
\end{equation*}

\noindent If the angle between $u$ and some root of $\mathfrak{q}_1$ is $45^\circ$, then by (\textbf{2.9}) and above
\begin{equation*}
2\geq ||u_1||^2+||u_2||^2\geq \frac{2n}{n+1}+\frac{n_2}{n_2+1}\geq \frac{2n}{n+1}+\frac{1}{2},
\end{equation*}
which is impossible when $n\geq 4$.\\

\noindent If the angle between $u$ and some root of $\mathfrak{q}_1$ is $60^\circ$ and the angle between $u$ and some root of $\mathfrak{q}_2$ is $45^\circ$, then by the calculations above,
\begin{equation*}
2\geq ||u_1||^2+||u_2||^2\geq ||u_1||^2+\frac{2n_2}{n_2+1}\geq ||u_1||^2+1.
\end{equation*}
\noindent By (\textbf{2.7}) and $n\geq 4$, the only length of weight that is less than or equal to $1$ is $\sqrt{\frac{n}{n+1}}$. Therefore, $||u_1||^2=\frac{n}{n+1}$. But then\\
\begin{equation*}
2\geq ||u_1||^2+||u_2||^2\geq \frac{n}{n+1}+\frac{2n_2}{n_2+1}\geq \frac{n}{n+1}+1.
\end{equation*}
\noindent This implies 
\begin{equation*}
||u||^2-(||u_1||^2+||u_2||^2)\leq 2-(\frac{n}{n+1}+1)=\frac{1}{n+1}.
\end{equation*}
\noindent Note that $\frac{1}{n+1}<\frac{1}{2}$ when $n\geq 2$. Therefore, by the last line of $(\textbf{2.9})$, we conclude that 
\begin{equation*}
||u||^2=||u_1||^2+||u_2||^2,
\end{equation*}
and hence, 
\begin{equation*}
2=\frac{n}{n+1}+\frac{2(n_2-k_2)(k_2+1)}{n_2+1}.
\end{equation*}
It is easy to check that this is not solvable for $n\geq 1$. Therefore, there does not exist a root $v\in\Phi_\mathfrak{g}$ such that the angle between $v$ and $u$ is $45^\circ$. Using similar arguments, we see if $n$ is big (at least $\geq 4$), then $u$ can only have two components. i.e.
\begin{equation*}
2=||u||^2=||u_1||^2+||u_2||^2.
\end{equation*} 
\noindent Since no weight of $A_n$ has length $1$ for $n\geq 4$, one of the $||u_1||^2$, $||u_2||^2$ is smaller than $1$ and the other is bigger than $1$. Since $\sqrt{\frac{m}{m+1}}$ is the only length of weight smaller than $1$ in $A_m$, together with (\textbf{2.9}), this gives the following Diophantine equation for positive integers $m,l,k$:\\
\begin{equation}
2=\frac{m}{m+1}+\frac{(l-k)(k+1)}{l+1}.\tag{D}
\end{equation}
~\bigskip

\noindent\textbf{Proposition 2.11.} The solutions of Diophantine equation (D) are $(m,l,k)=(1,5,2)$, $(1,7,1)$, $(1,7,5)$, $(2,5,1)$, $(2,5,3)$, $(4,4,1)$ and $(4,4,2)$.\\

\noindent \textbf{Proof.} Write $a=l-k$ and $b=k+1$. Since $m$ is positive, thus we have 
\begin{equation*}
\frac{ab}{a+b}<2
\end{equation*}
which is equivalent to 
\begin{equation*}
(a-2)(b-2)\leq 3.
\end{equation*}
If $a-2$ and $b-2$ are non-zero, then there are only finitely many possible pairs $(a,b)$. In this case, there exists only one solution $(m,l,k)=(1,5,2)$ which corresponds to $(a,b)=(3,3)$.\\

\noindent If $a-2=0$, then the equation becomes
\begin{equation*}
2=\frac{m}{m+1}+\frac{2b}{b+2}
\end{equation*}
which is equivalent to
\begin{equation*}
m(b-2)=4.
\end{equation*}
Therefore, we see $(m,b)=(1,6),(2,4),(4,3)$ and we deduce easily that $(m,l,k)=(1,7,1)$, $(1,7,5)$, $(2,5,1)$, $(2,5,3)$, $(4,4,1)$, $(4,4,2)$.\qed\\

\noindent(\textbf{2.12}) From Proposition 2.11, we see that the $m$, $l$ coordinates of solutions of (D) do not belong to $\{6,8,9,10,11,...\}$. Therefore, there does not exist a root $u\in\Phi_\mathfrak{h}$ such that $u\notin (\Lambda_\mathfrak{q}\otimes\mathbb{R})\cup (\Lambda_\mathfrak{q}\otimes\mathbb{R})^\bot$ for any $A_n$ factor $\mathfrak{q}$ of $\mathfrak{g}$ where $n\in\{6,8,9,10,11,...\}$. In other words, any root $u\in\Phi_\mathfrak{h}$ must either be contained in or perpendicular to $\Lambda_\mathfrak{q}\otimes\mathbb{R}\subset \Lambda_\mathfrak{g}\otimes\mathbb{R}$. On the other hand, since any simple root system is irreducible, there exist simple factors $\mathfrak{p}_1,...,\mathfrak{p}_s$ of $\mathfrak{h}$ such that $\Lambda_\mathfrak{q}\otimes\mathbb{R}=(\Lambda_{\mathfrak{p}_1}\otimes\mathbb{R})\oplus\cdots\oplus (\Lambda_{\mathfrak{p}_s}\otimes\mathbb{R})$. If we consider the sum of all the weights in the formal character sitting in this subspace of $t_\mathfrak{g}^*$, then it is actually a formal character coming from some faithful representations of $\mathfrak{q}$ and $\mathfrak{p}_1\oplus\cdots\oplus\mathfrak{p}_s$ \cite[Exercise 23.42]{FH}. Thus, we reduce our problem to the case that $\mathfrak{g}$ is a simple complex Lie algebra of type $A_n$, $n\in\{6,8,9,10,11,...\}$. If we can prove that $s=1$, then $\mathfrak{p}_1=\mathfrak{q}$ and the number of $A_n$ factors of $\mathfrak{g}$ and $\mathfrak{h}$ for $n\in\{6,8,9,10,11,...\}$ are the same. \\

\noindent(\textbf{2.13}) Suppose $\mathfrak{g}$ is simple of type $A_n$, $n\in\{4,5,6,9,10,11,...\}$ and $\mathfrak{h}=\mathfrak{p}_1\oplus\cdots\oplus\mathfrak{p}_s$ is a direct sum of simple factors of type $A_n$. Assume the length of roots of $\mathfrak{g}$ is $\sqrt{2}$. Let $u$ be a root of $\mathfrak{h}$ and $\theta$ the angle between $u$ and some root $v\in\Phi_\mathfrak{g}$. Consider the following $3$ cases:
 \begin{enumerate}
 \item [$\theta=30^\circ$] This is impossible because we can always choose some root $v'$ (since $n\geq 3$) of $\mathfrak{g}$ that makes an angle of $30^\circ$ with $u$ such that $\{u,v,v'\}$ spans a $3$-dimensional space and then use the argument in Lemma 2.6. \\
 
 \item [$\theta=45^\circ$] Then $u$ is either of length $1$ or $2$ by (\textbf{2.4}). It is impossible for $u$ and some root of $\mathfrak{g}$ to be parallel because $\theta=45^\circ$. Hence, the angle between $u$ and any root of $\mathfrak{g}$ belongs to $\{45^\circ,90^\circ,135^\circ\}$. By Proposition 2.8 and computations in (\textbf{2.9}), we deduce that

\begin{equation*}
4=||u||^2=\frac{4(n-k)(k+1)}{n+1}\hspace{.3in}\mathrm{or}\hspace{.3in}1=||u||^2=\frac{(n-k)(k+1)}{n+1}.
\end{equation*}
Hence, we for both cases obtain the equation
\begin{equation*}
k(n-k-1)=1
\end{equation*}
which implies $(n,k)=(3,1)$. This is impossible since $3\notin \{4,5,6,9,10,11,...\}$.  \\

\item [$\theta=60^\circ$] Length of $u$ is $\sqrt{2}$ by (\textbf{2.4}). If $u$ is not equal to any root of $\mathfrak{g}$, then the angle between $u$ and any root of $\mathfrak{g}$ belongs to $\{60^\circ,90^\circ,120^\circ\}$. By Proposition 2.8 and computations in (\textbf{2.9}),
\begin{equation*}
2=||u||^2=\frac{(n-k)(k+1)}{n+1}.
\end{equation*}
\begin{equation*}
2(n+1)=(n-k)(k+1).
\end{equation*}
Let $N=n+1$ and $K=k+1$, we have
\begin{equation*}
K^2-NK+2N=0.
\end{equation*}
Since $N^2-8N$, the discriminant of the left hand side is a perfect square $\Delta^2$, thus we have
\begin{equation*}
(N-\Delta)(N+\Delta)=8N.
\end{equation*}

\noindent Assume $\Delta$ is non-negative, we have $1\leq (N-\Delta)\leq 8$. So we only have $8$ cases to consider: $(N-\Delta)\in\{1,2,...,8\}$. We find that $N$ can only be $8,9$, thus $n$ can only be $7,8$ which do not belong to $\{4,5,6,9,10,...\}$.\end{enumerate}

\noindent So $\theta=0^\circ$ or $180^\circ$ and we conclude that any root $u$ of $\mathfrak{h}$ lies on the line spanned by a root $\mathfrak{g}$. If we can prove that $\mathfrak{h}$ is simple, then $\mathfrak{g}=\mathfrak{h}$ and we are done. Write $\mathfrak{h}=\mathfrak{p}_1\oplus\mathfrak{p}_2\oplus\cdots\oplus\mathfrak{p}_s$ and let $n_i$ be the rank of $\mathfrak{p}_i$. If we take a base $S_i$ of $\mathfrak{p}_i$ for all $1\leq i\leq s$, then $S_i$ is a subset of $\Phi_\mathfrak{g}$ for all $i$ and $S_1$ is orthogonal to $S_i$ for $i\geq 2$. Note that the union of all $S_i$ forms a basis of $\Lambda_\mathfrak{g}\otimes\mathbb{R}$, therefore $S_2\cup\cdots\cup S_s$ spans a space of dimension $n-n_1$. By (\textbf{2.7}), the set of roots of $\mathfrak{g}$ is given by $\{e_i-e_j:1\leq i\neq j\leq n+1\}$. Without loss of generality, we may assume 
\begin{equation*}
S_1=\{c(e_1-e_2),c(e_2-e_3),...,c(e_{n_1}-e_{n_1+1})\}
\end{equation*}
for some constant $c\in\{1/2,1,2\}$ by (\textbf{2.4}). If $a_1e_1+\cdots+a_ne_n$ is orthogonal to $S_1$, then $a_1=a_2=\cdots=a_{n_1}=a_{n_1+1}$ by (\textbf{2.7}). Therefore, the set of roots of $\mathfrak{g}$ that are orthogonal to $S_1$ is 
\begin{equation*}
\{e_i-e_j:~n_1+1\leq i\neq j\leq n+1\}
\end{equation*}
which spans a space of dimension $n-n_1-1$. It contradicts that $S_2\cup\cdots\cup S_s$ spans a space of dimension $n-n_1$. So we conclude that $\mathfrak{h}$ is simple and is equal to $\mathfrak{g}=A_n$ if $n\in\{4,5,6,9,10,11,...\}$. If we consider the intersection of the two sets $\{6,8,9,10,11,...\}$ and $\{4,5,6,9,10,11,...\}$ in (\textbf{2.12}) and (\textbf{2.13}), then we can get the following theorem.\\

\noindent \textbf{Theorem 2.14} If faithful representations of two complex semisimple Lie algebras $\mathfrak{g}$ and $\mathfrak{h}$ have the same formal character, then the number of $A_n$ factors of $\mathfrak{g}$ and $\mathfrak{h}$ are the same when $n\in\{6,9,10,...\}$.\\

\noindent \textbf{Proof.} First, we may assume each faithful representation contains the adjoint representation as a subrepresentation by (\textbf{2.2}). Hence the roots $\Phi_\mathfrak{g}$, $\Phi_\mathfrak{h}$ really appear in the formal character. Since we only care about the number of $A_n$ factors for $n\in\{6,9,10,...\}$, we assume $\mathfrak{g}$ and $\mathfrak{h}$ only have simple factors of type $A_n$ by Lemma 2.5. Now if $\mathfrak{q}$ is a simple factor of $\mathfrak{g}$ with rank $n\in\{6,9,10,...\}$ and $u\in\Phi_\mathfrak{h}$, (\textbf{2.12}) says that $u$ is either contained in or perpendicular to $\Lambda_\mathfrak{q}\otimes\mathbb{R}$. So we may reduce to the case that $\mathfrak{g}$ is $\mathfrak{q}$ and if we can prove that $\mathfrak{h}$ is also simple, then we are done. Finally (\textbf{2.13}) finishes that part.\qed\\

\noindent (\textbf{2.15}) We are going to prove that the parities of the number of $A_4$ factors of $\mathfrak{g}$ and $\mathfrak{h}$ are equal. Again we may assume our Lie algebras only have simple factors of type $A_n$ by Lemma 2.5 and the roots of $\mathfrak{h}$ appear in the formal character by (\textbf{2.2}). Therefore, the roots of $\mathfrak{h}$ belong to $\Lambda_\mathfrak{g}$. Suppose 
\[ \begin{array}{lcl}

\mathfrak{g}=\mathfrak{g}_1\oplus\mathfrak{g}_2\\
\mathfrak{h}=\mathfrak{h}_1\oplus\mathfrak{h}_2\\
\end{array}\]
are decompositions into a direct sum of proper ideals such that $\Lambda_{\mathfrak{g}_1}\otimes\mathbb{R}=\Lambda_{\mathfrak{h}_1}\otimes\mathbb{R}$. Then we also have $\Lambda_{\mathfrak{g}_2}\otimes\mathbb{R}=\Lambda_{\mathfrak{h}_2}\otimes\mathbb{R}$ by orthogonality. If this is the case, then we say the pair $(\mathfrak{g}\curvearrowright V,\mathfrak{h}\curvearrowright V')$ of representations is \emph{reducible}. It is because the weights of the original formal character sitting inside $\Lambda_{\mathfrak{g}_i}\otimes\mathbb{R}=\Lambda_{\mathfrak{h}_i}\otimes\mathbb{R}$ come from some faithful representations of $\mathfrak{g}_i$ and $\mathfrak{h}_i$, $i=1,2$. We can reduce our problem to some pairs $(\mathfrak{g}_i\curvearrowright V_i,\mathfrak{h}_i\curvearrowright V_i')$ for $i=1,2$.  Now we assume that $(\mathfrak{g}\curvearrowright V,\mathfrak{h}\curvearrowright V')$ is not reducible. If $v\in\Phi_\mathfrak{g}$ and $u\in\Phi_\mathfrak{h}$, then by (\textbf{2.4}) we have 
\begin{equation*}
\frac{\left\langle v, v \right\rangle}{\left\langle u, u \right\rangle}\in 2^a3^b
\end{equation*}
\noindent for some integers $a,b$.\\

\noindent (\textbf{2.16}) Assume the pair $(\mathfrak{g}\curvearrowright V,\mathfrak{h}\curvearrowright V')$ is not reducible. Let $n$ be their common rank and write
\[ \begin{array}{lcl}

\mathfrak{g}=\mathfrak{q}_1\oplus\cdots\oplus\mathfrak{q}_m\\
\mathfrak{h}=\mathfrak{p}_1\oplus\cdots\oplus\mathfrak{p}_s\\
\end{array}\]
\noindent as direct sums of simple factors. Let the rank of $\mathfrak{q}_i$, $\mathfrak{p}_j$ be $n_i$, $r_j$ respectively. The standard basis for the weight lattice of $\mathfrak{q}_i$, $\Lambda_{\mathfrak{q}_i}$, $1\leq i\leq m$ is given by
\begin{equation*}
\{e^i_1,e^i_2,....,e^i_{n_i}\}.
\end{equation*}
\noindent So a basis for the weight lattice of $\mathfrak{g}$, $\Lambda_\mathfrak{g}$ is given by
\begin{equation*}
B=\{e^1_1,e^1_2,....,e^1_{n_1},e^2_1,e^2_2,....,e^2_{n_2},....,e^m_1,e^m_2,....,e^m_{n_m}\}.
\end{equation*}
\noindent If we normalize the inner product $\left\langle ~, ~ \right\rangle$ so that the length of roots of $\mathfrak{q}_1$ is $\sqrt{2}$, then by (\textbf{2.15}) the positive definite $n\times n$ matrix $Q$ defining $\left\langle ~, ~ \right\rangle$ on $\Lambda_{\mathfrak{g}}\otimes\mathbb{R}$ under the basis $B$ is 

\begin{equation*}
   Q=\begin{pmatrix}
      Q_1 &0 &... &0\\
      0 &\gamma_2Q_2  &... &0\\
      \vdots &\vdots &\ddots &\vdots\\
      0 &0 &... &\gamma_mQ_m
\end{pmatrix},
\end{equation*}
\noindent where $\gamma_i\in 2^\mathbb{Z}3^\mathbb{Z}$ by (\textbf{2.15}) and $Q_i$ is an $n_i\times n_i$ matrix in the following form by (\textbf{2.7}) for $1\leq i\leq m$:
\begin{equation*}
   Q_i=\begin{pmatrix}
      \frac{n_i}{n_i+1} &\frac{-1}{n_i+1} &... &\frac{-1}{n_i+1}\\
      \frac{-1}{n_i+1} &\frac{n_i}{n_i+1}  &... &\frac{-1}{n_i+1}\\
      \vdots &\vdots &\ddots &\vdots\\
      \frac{-1}{n_i+1} &\frac{-1}{n_i+1} &... &\frac{n_i}{n_i+1}
\end{pmatrix}.
\end{equation*}
\noindent If we take a base $S=\{u_1,...,u_n\}$ of root system of $\mathfrak{h}$, then it defines an $n\times n$ matrix $A$ whose $i$-th column is $u_i$ in terms of the basis $B$. Since all $u_i$ belong to the weight lattice of $\mathfrak{g}$, $A$ is a matrix of integral entries. If we enumerate $S$ in a proper way, then we have the following equation
\begin{equation*}
   A^tQA=\begin{pmatrix}
      \mu_1C_1 &0 &... &0\\
      0 &\mu_2C_2  &... &0\\
      \vdots &\vdots &\ddots &\vdots\\
      0 &0 &... &\mu_sC_s
\end{pmatrix},
\end{equation*}
\noindent where $\mu_j\in 2^\mathbb{Z}3^\mathbb{Z}$ by (\textbf{2.15}) and $C_j$ is the Cartan matrix of $\mathfrak{p_j}$ for $1\leq j\leq s$. We know that the determinant of the Cartan matrix $C_j$ is $r_j+1$ and it is easy to check that the determinant of $Q_i$ is $\frac{1}{n_i+1}$. Therefore, by taking determinant of the matrix equation, we get
\begin{equation*}
\frac{\gamma_2....\gamma_m \mathrm{det}(A)^2}{(n_1+1)....(n_m+1)}=\mu_1....\mu_s(r_1+1)....(r_s+1).
\end{equation*}
\noindent We see det$(A)$ is an integer and $\gamma_i,\mu_j\in 2^\mathbb{Z}3^\mathbb{Z}$. If we reduce this equation modulo $(\mathbb{Q}^*)^2$, then by the unique factorization of primes we see the parities of $\mathrm{ord}_5(n_1+1)....(n_m+1)$ and  $\mathrm{ord}_5(r_1+1)....(r_s+1)$ are the same. By Theorem 2.14, the number of $A_n$ factors for $n\geq 9$ is invariant, so this implies the parity of number of $A_4$ factors is an invariant because $4+1=5$. Thus, we have just proved the following theorem.\\

\noindent \textbf{Theorem 2.17} If faithful representations of two complex semisimple Lie algebras $\mathfrak{g}$ and $\mathfrak{h}$ have the same formal character, then the parities of the numbers of $A_4$ factors of $\mathfrak{g}$ and $\mathfrak{h}$ are the same.\\

\noindent \textbf{Definition 2.18} Consider the free abelian group $\mathfrak{F}$ of virtual complex simple Lie algebras. Semisimple Lie algebras are subset $\mathfrak{F}$ in a natural way. We divide by the subgroup $\mathfrak{D}$ generated by all expressions $\mathfrak{g}-\mathfrak{h}$, where $\mathfrak{h}\subset\mathfrak{g}$ are semisimple of same rank. We say that two complex semisimple Lie algebras $\mathfrak{g}$ and $\mathfrak{h}$ satisfy \emph{equal-rank subalgebra equivalence}, denoted by $\mathfrak{g}\approx\mathfrak{h}$, if they have same image in $\mathfrak{F}/\mathfrak{D}$.\\

\noindent Theorem 2.19 follows easily from the results above.\\

\noindent \textbf{Theorem 2.19} If faithful representations of two complex semisimple Lie algebras $\mathfrak{g}$ and $\mathfrak{h}$ have the same formal character, then $\mathfrak{g}$ and $\mathfrak{h}$ satisfy equal-rank subalgebra equivalence.\\

\noindent \textbf{Proof.} By Lemma 2.5 we may assume our Lie alegbras consist of simple factors of type $A_n$. By Theorem 2.14, we further assume the simple factors are of types $A_n$, $n\in\{1,2,3,4,5,7,8\}$. From \textbf{Table $1$}, the following Lie algebras are equivalent:
\begin{enumerate}
\item $A_4\times A_4\approx A_8\approx A_2\times A_2\times A_2\times A_2$;
\item $A_7\approx A_2\times A_5$;
\item $A_1\times A_5\approx A_2\times A_2\times A_2$;
\item $A_2\approx A_1\times A_1$.
\end{enumerate}
\noindent Hence by Theorem 2.17 and the list above, we may assume $n\in\{1,3\}$ and it suffices to prove that $A_3\approx A_1\times A_1\times A_1$. But we know that
\begin{equation*}
A_3\times A_1\approx \mathfrak{so}(6)\times\mathfrak{so}(3)\approx\mathfrak{so}(4)\times\mathfrak{so}(5)\approx\mathfrak{so}(4)\times\mathfrak{so}(4)\approx A_1\times A_1\times A_1\times A_1.
\end{equation*}
\noindent So we are done.\qed\\

\noindent \textbf{Proposition 2.20} Two semisimple Lie algebras $\mathfrak{g}$ and $\mathfrak{h}$ satisfy equal-rank subalgebra equivalence if and only if the number of $A_n$ factors when $n\in\{6,9,10,11,...\}$ and the parity of the number of $A_4$ factors of them are the same.

\noindent \textbf{Proof.} We have seen that equal-rank subalgebra equivalence is implied by the conditions on $A_n$ factors in the proof of the above theorem. Suppose two complex semisimple Lie algebras $\mathfrak{g}$ and $\mathfrak{h}$ of the same rank satisfy equal-rank subalgebra equivalence, then in the free group $\mathfrak{F}$ we have
\begin{equation*}
\mathfrak{g}-\mathfrak{h}=\mathfrak{h}_1-\mathfrak{g}_1+\mathfrak{h}_2-\mathfrak{g}_2+\cdots+\mathfrak{h}_k-\mathfrak{g}_k
\end{equation*}
where $\mathfrak{h}_i\subset\mathfrak{g}_i$ or $\mathfrak{g}_i\subset\mathfrak{h}_i$ and are of the same rank for all $i$. Thus we obtain
 \begin{equation*}
\mathfrak{g}+\mathfrak{g}_1+\cdots+ \mathfrak{g}_k=     \mathfrak{h}+\mathfrak{h}_1+\cdots+\mathfrak{h}_k
\end{equation*}
By Theorem 2.14 and 2.17, one sees easily by taking a faithful representation of $\mathfrak{g}_i$ or $\mathfrak{h}_i$ that for each $i$, $\mathfrak{g}_i$ and $\mathfrak{h}_i$ have the the number of $A_n$ factors for $n\in\{6,9,10,...\}$ and the parity of the numbers of $A_4$ factors. So do $\mathfrak{g}$ and $\mathfrak{h}$ by the equation above.\qed\\

\textbf{Table $1$} \cite[Table 5]{GOV}. The table lists some semisimple maximal subalgebras $\mathfrak{f}$ of maximal rank in simple complex Lie algebra $\mathfrak{g}$ (up to conjugacy in $\mathfrak{g}$).
\begin{center}
\begin{tabular}{|c|c|} \hline
\hspace{.5in}$\mathfrak{g}$\hspace{.5in}  & \hspace{.9in}$\mathfrak{f}$\hspace{.9in} \\ \hline
$\mathfrak{so}_{2l+1}(\mathbb{C})$  & $\mathfrak{so}_{2k}(\mathbb{C})\oplus\mathfrak{so}_{2(l-k)+1}(\mathbb{C})$ \\ 
$l\geq2$  & $2\leq k\leq l$ \\ \hline
$\mathfrak{sp}_{2l}(\mathbb{C})$ & $\mathfrak{sp}_{2k}(\mathbb{C})\oplus\mathfrak{sp}_{2(l-k)}(\mathbb{C})$ \\ 
$l\geq3$  & $1\leq k\leq [\frac{l}{2}]$ \\ \hline
$\mathfrak{so}_{2l}(\mathbb{C})$ & $\mathfrak{so}_{2k}(\mathbb{C})\oplus\mathfrak{so}_{2(l-k)}(\mathbb{C})$ \\ 
$l\geq4$  & $2\leq k\leq \frac{l+1}{2}$ \\ \hline
$E_6$  & $\mathfrak{sl}_2(\mathbb{C})\oplus\mathfrak{sl}_6(\mathbb{C})$ \\ 
 ~ & $\mathfrak{sl}_3(\mathbb{C})\oplus\mathfrak{sl}_3(\mathbb{C})\oplus\mathfrak{sl}_3(\mathbb{C})$ \\ \hline
  & $\mathfrak{sl}_2(\mathbb{C})\oplus\mathfrak{so}_{12}(\mathbb{C})$ \\ 
$E_7$ & $\mathfrak{sl}_3(\mathbb{C})\oplus\mathfrak{sl}_{6}(\mathbb{C})$\\ 
& $\mathfrak{sl}_8(\mathbb{C})$\\ \hline

 & $\mathfrak{sl}_2(\mathbb{C})\oplus E_7$\\ 
 & $\mathfrak{sl}_3(\mathbb{C})\oplus E_6$\\ 
$E_8$ & $\mathfrak{sl}_5(\mathbb{C})\oplus\mathfrak{sl}_{5}(\mathbb{C})$\\ 
 & $\mathfrak{so}_{16}(\mathbb{C})$\\ 
  & $\mathfrak{sl}_9(\mathbb{C})$ \\ \hline

  & $\mathfrak{sl}_2(\mathbb{C})\oplus \mathfrak{sp}_6(\mathbb{C})$ \\
$F_4$  & $\mathfrak{sl}_3(\mathbb{C})\oplus\mathfrak{sl}_{3}(\mathbb{C})$ \\ 
  & $\mathfrak{so}_9(\mathbb{C})$ \\ \hline

$G_2$  & $\mathfrak{sl}_3(\mathbb{C})$ \\ 
  & $\mathfrak{so}_4(\mathbb{C})$ \\ \hline
\end{tabular}
\end{center}
~\\\\

\begin{center}
\textbf{$\mathsection3$ Compatible system of $\ell$-adic representations and $\ell$-independence}
\end{center}

\noindent(\textbf{3.1}) We follow the terminology and notations of Serre \cite{Serre3}. Let $K$ be a number field, $G_K$ its absolute Galois group, $\sum_K$ the set of all finite places of $K$ and $\ell$ a prime number. An \emph{$\ell$-adic representation of $K$} is a continuous homomorphism $\rho: G_K\rightarrow \mathrm{GL}_n(\mathbb{Q}_\ell)$ for some $n$. Here $G_K$ is equipped with the profinite topology and $\mathrm{GL}_n(\mathbb{Q}_\ell)$ is an $\ell$-adic Lie group. An \emph{abelian $\ell$-adic representation of $K$} is an $\ell$-adic representation of $G_K$ that factors through $G_K^\mathrm{ab}$.\\

\noindent(\textbf{3.2}) Let $v\in \sum_K$ and $p_v$ the characteristic of the residue field $k_v$ of the place $v$. If $w$ is a valuation of $\overline{K}$ extending $v$, we denote the decomposition group, inertia group and Frobenius element of $w$ by $D_w$, $I_w$ and $F_w$ respectively. $D_w$ and $I_w$ are closed subgroups of $G_K$. We say that $\rho$ is unramified at $v$ if $\rho(I_w)$ is trivial for any valuation $w$ of $\overline{K}$ extending $v$. If the representation $\rho$ is unramified at $v$, then the restriction of $\rho$ to $D_w$ factors through $D_w/I_w$ and $\rho(F_w)$ is defined for any $w|v$ . We denote it by $F_{w,\rho}$ and the conjugacy class of $F_{w,\rho}$ in $\mathrm{GL}_n(\mathbb{Q}_\ell)$ by $F_{v,\rho}$.\\

\noindent\textbf{Definition 3.3.} An $\ell$-adic representation $\rho$ is said to be \emph{rational} if there exists a finite subset $S$ of $\sum_K$ such that:
\begin{enumerate}
\item[(a)] Any element of $\sum_K\setminus S$ is unramified with respect to $\rho$.
\item[(b)] If $v\notin S$, the coefficients of $P_{v,\rho}(T):=\mathrm{det}(1-F_{v,\rho}T)$ belongs to $\mathbb{Q}$.
\end{enumerate}

\noindent\textbf{Definition 3.4.} Let $\ell'$ be a prime, $\rho'$ an $\ell'$-adic representation of $K$, and assume that $\rho$, $\rho'$ are rational. Then $\rho$, $\rho'$ are said to be \emph{compatible} if there exists a finite subset $S$ of $\sum_K$ such that $\rho$ and $\rho'$ are unramified outside of $S$ and $P_{v,\rho}(T)=P_{v,\rho'}(T)$ for $v\in\sum_K\setminus S$.\\

\noindent\textbf{Definition 3.5.} Let $\mathscr{P}$ be the set of prime numbers. For each prime $\ell$ let $\rho_\ell$ be a rational $\ell$-adic representation of $K$. The system $\{\rho_\ell\}_{\ell\in \mathscr{P}}$ is said to be \emph{compatible} if $\rho_\ell$, $\rho_{\ell'}$ are compatible for any two primes $\ell$, $\ell'$. The system $\{\rho_\ell\}_{\ell\in \mathscr{P}}$ is said to be \emph{strictly compatible} if there exists a finite subset $S$ of $\sum_K$ such that:
\begin{enumerate}
\item[(a)] Let $S_\ell=\{v:~p_v=\ell\}$. Then, for every $v\notin S\cup S_\ell$, $\rho_\ell$ is unramified at $v$ and $P_{v,\rho_\ell}(T)$ has rational coefficients.
\item[(b)] $P_{v,\rho_\ell}(T)=P_{v,\rho_{\ell'}}(T)$ if $v\notin S\cup S_\ell\cup S_{\ell'}$.\\

\end{enumerate}

\noindent(\textbf{3.6}) Serre associates to every number field $K$ a projective family $\{S_\mathfrak{m}\}$ of commutative algebraic groups over $\mathbb{Q}$ and shows that each $S_\mathfrak{m}$ gives rise to a compatible system of rational $\ell$-adic representation of $K$ (\cite[Chap. 2, $\mathsection$1, 2]{Serre3}). We give a brief introduction to the algebraic groups $S_\mathfrak{m}$ associated to the number field $K$. Let $S$ be a finite subset of $\sum_K$. Then by a \emph{modulus of support} $S$ we mean a family $\mathfrak{m}=(m_v)_{v\in S}$ where the $m_v$ are integers $\geq 1$. If $E$ is the group of units of $K$, 
\begin{equation*}
E_\mathfrak{m}:=\{u\in E: v(1-u)\geq m_v~for~all~v\in S\}.
\end{equation*}
Then $S_\mathfrak{m}$ is an algebraic group over $\mathbb{Q}$ whose connected component is 
\begin{equation*}
T_\mathfrak{m}:=\mathrm{Res}_{K/\mathbb{Q}}(\mathbb{G}_{m/K})/\overline{E_\mathfrak{m}},
\end{equation*}

\noindent where $\mathrm{Res}_{K/\mathbb{Q}}(\mathbb{G}_{m/K})$ is obtained from the multiplicative group $\mathbb{G}_m$ by restriction of scalars from $K$ to $\mathbb{Q}$ and $\overline{E_\mathfrak{m}}$ is the Zariski closure of $E_\mathfrak{m}$ in the algebraic torus $\mathrm{Res}_{K/\mathbb{Q}}(\mathbb{G}_{m/K})$. $S_\mathfrak{m}$ are called the \emph{Serre groups} associated to $K$. \\

\noindent\textbf{Definition 3.7.} Let $H$ be a linear algebraic group over $\mathbb{Q}$, and let $K$ be a number field. A continuous homomorphism $\rho:G_K\rightarrow H(\mathbb{Q}_\ell)$ is called an \emph{$\ell$-adic representation of $K$ with values in $H$}.\\

 \noindent One defines in an analogous way what it means for $\rho$ to be unramified at a place or rational and for a system $\{\rho_\ell\}$ to be compatible or strictly compatible. Using class field theory, Serre defines for each prime $\ell$ and for each modulus $\mathfrak{m}$ an $\ell$-adic representations $\epsilon_\ell$ with values in $S_\mathfrak{m}$ (Definition 3.7)
\begin{equation*}
\epsilon_\ell: G_K^\mathrm{ab}\longrightarrow S_\mathfrak{m}(\mathbb{Q}_\ell).
\end{equation*}
 
\noindent \textbf{Theorem 3.8.}\begin{enumerate}
\item The dimensions of the Serre groups associated to $K$ only depend on $K$. Denote the common dimension by $d_K$.
\item \cite[Chap. 2.2.3]{Serre3} The image of $\epsilon_\ell$ is an open subgroup of $S_\mathfrak{m}(\mathbb{Q}_\ell)$ and is Zariski dense in $S_\mathfrak{m}$.
\item \cite[Chap. 2.2.5]{Serre3} Let $\phi:S_\mathfrak{m}\rightarrow \mathrm{GL}_{n,\mathbb{Q}}$ be a $\mathbb{Q}$-morphism, then the representation
\begin{equation*}
\phi\circ\epsilon_\ell:G_K^\mathrm{ab}\rightarrow \mathrm{GL}_n(\mathbb{Q}_\ell)
\end{equation*}
is semisimple and $\{\epsilon_\ell\}_{\ell\in \mathscr{P}}$ is strictly compatible.
\item The dimension of $S_\mathfrak{m}(\mathbb{Q}_\ell)$ as an $\ell$-adic Lie group is equal to $d_K$ for all $\ell$.
\end{enumerate}

\noindent \textbf{Proof.} We have provided references for (2) and (3), so we just need to prove (1) and (4). Since $E_\mathfrak{m}$ is of finite index in $E$ \cite[Chap. 2.2.1]{Serre3}, Serre groups have the same dimension. By taking an $\mathbb{Q}$-embedding $\phi$ of $S_\mathfrak{m}$ to some $\mathrm{GL}_{n,\mathbb{Q}}$, one obtains (4) by the algebraicity of the Lie algebra of $\phi\circ\epsilon_\ell(G_K)$ \cite[Thm. 4]{Hen}, (2) and (3) of this theorem.\qed\\

\noindent(\textbf{3.9}) We are going to define local algebraicity for an abelian $\ell$-adic representation of a number field $K$. We first need to define local algebraicity for $K$, a finite extension of $\mathbb{Q}_\ell$. Let $T=\mathrm{Res}_{K/\mathbb{Q}_\ell}(\mathbb{G}_{m/K})$ be the algebraic torus over $\mathbb{Q}_\ell$ and $V$ a finite dimensional vector space over $\mathbb{Q}_\ell$. If $i:K^*\rightarrow G_K^\mathrm{ab}$ is the canonical homomorphism of local class field theory and $\rho: G_K^\mathrm{ab}\rightarrow \mathrm{Aut}(V)$ is an abelian $\ell$-adic representation of $K$, we then get a continuous homomorphism $\rho\circ i$ of $K^*=T(\mathbb{Q}_\ell)$ into $\mathrm{Aut}(V)$. \\

\noindent\textbf{Definition 3.10.} The representation $\rho$ is said to be \emph{locally algebraic} if there is an algebraic morphism $r: T\rightarrow \mathrm{GL}_V$ such that $\rho\circ i(x)=r(x^{-1})$ for all $x\in K^*$ close enough to $1$.\\

\noindent(\textbf{3.11}) Now if $K$ is a number field and $V_\ell$ is a finite dimensional vector space over $\mathbb{Q}_\ell$. Let 
\begin{equation*}
\rho: G_K^\mathrm{ab}\longrightarrow \mathrm{Aut}(V_\ell)
\end{equation*}

\noindent be an abelian $\ell$-adic representation of $K$. Let $v\in \sum_K$ be a place of $K$ of residue characteristic $\ell$ and let $D_v\subset G_K^\mathrm{ab}$ be the corresponding decomposition group which is isomorphic to $\mathrm{Gal}(\overline{K_v}/K_v)^\mathrm{ab}$. Hence, we get an $\ell$-adic representation of $K_v$ by composition
\begin{equation*}
\rho_v: \mathrm{Gal}(\overline{K_v}/K_v)^\mathrm{ab}\rightarrow D_v\rightarrow \mathrm{Aut}(V_\ell).
\end{equation*} 

\noindent\textbf{Definition 3.12.} The representation $\rho$ is said to be \emph{locally algebraic} if all the local representations $\rho_v$, with $p_v=\ell$, are locally algebraic in the sense of Definition 3.10.\\

\noindent The following two theorems are crucial for this section. Let $\rho: G_K^\mathrm{ab}\rightarrow \mathrm{Aut}(V_\ell)$ be an abelian $\ell$-adic representation of the number field $K$.\\

\noindent\textbf{Theorem 3.13.} \cite[Chap. 3, Thm. 2]{Serre3} If $\rho$ is locally algebraic, then there exists a modulus $\mathfrak{m}$, an abelian $\ell$-adic representation $\epsilon_\ell: G_K^\mathrm{ab}\rightarrow S_\mathfrak{m}(\mathbb{Q}_\ell)$ (\textbf{3.8}) and a morphism of algebraic groups $\phi: S_\mathfrak{m}\times_\mathbb{Q}\mathbb{Q}_\ell\rightarrow \mathrm{GL}_{V_\ell}$ over $\mathbb{Q}_\ell$  such that $\rho=\phi\circ \epsilon_\ell$.\\

\noindent\textbf{Theorem 3.14.} \cite[Chap. 3 $\mathsection$3]{Serre3}, \cite{W}, \cite[$\mathsection$5]{Hen} If $\rho$ is rational and semisimple, then $\rho$ is locally algebraic.\\

\noindent\textbf{Remark 3.15.} An abelian $\ell$-adic representation $\rho$ of $K$ is always unramified at all $v\in\sum_K\setminus S$, for some finite set $S$ (\cite[Chap.3 $\mathsection$2.2]{Serre3}). The proof of Theorem 3.14 \cite[$\mathsection$5]{Hen} consists of two parts. Part one proves that $\rho$ is almost locally algebraic, i.e. there exists an integer $N$ such that $\rho^N$ is locally algebraic. Part two proves that if $\rho$ is almost locally algebraic, then $\rho$ is locally algebraic. The crucial observation is that part one does not use the full rational condition, it only needs $\rho(F_w)$'s eigenvalues are algebraic for all $w|v$ when $v\notin S$. We will use this in the proof of Proposition 3.18.\\

\noindent(\textbf{3.16}) Suppose $\rho:G_K\rightarrow \mathrm{GL}_n(\mathbb{Q}_\ell)$ is a semisimple, rational $\ell$-adic representation of $K$. Let $G_\ell$ be the Zariski closure of $\rho(G_K)$ in $\mathrm{GL}_{n,\mathbb{Q}_\ell}$. It is defined over $\mathbb{Q}_\ell$. Semisimplicity of $\rho$ implies $G_\ell$ is a reductive group. Indeed, $G_\ell$ acts on an $n$-dimensional vector space $V$ over $\mathbb{Q}_\ell$. Semisimplicity implies $V$ decomposes into a direct sum of irreducibles $V_1\oplus V_2\oplus\cdots\oplus V_m$. If $U$ is the unipotent radical of $G_\ell^\circ$, then the eigenspace $W$ (of eigenvalue $1$, the only eigenvalue) of $U$ decomposes as $W_1\oplus W_2\oplus\cdots\oplus W_m$ such that $W_i\subset V_i$ and $W_i$ is non-trivial for all $1\leq i\leq m$. Since $U$ is normal in $G_\ell$, $W$ is an invariant subspace of $V$. Therefore, $W=V$, and $G_\ell$ is a reductive algebraic group. \\

\noindent(\textbf{3.17}) We may now assume $G_\ell$ is connected, reductive which corresponds to restricting $\rho$ to the open subgroup $G_{K^\mathrm{conn}}$ of $G_K$. The quotient group $G_\ell/[G_\ell,G_\ell]$ is a $\mathbb{Q}_\ell$-torus. Let \begin{equation*}
j:G_\ell/[G_\ell,G_\ell]\rightarrow \mathrm{GL}_{m,\mathbb{Q}_\ell}
\end{equation*}
\noindent be an embedding defined over $\mathbb{Q}_\ell$ and introduce the map 
\begin{equation*}
\theta:G_\ell\rightarrow G_\ell/[G_\ell,G_\ell]\stackrel{j}{\rightarrow} \mathrm{GL}_{m,\mathbb{Q}_\ell}.
\end{equation*}
\noindent If $\rho$ is unramified outside a finite subset $S$ of $\sum_K$ and $w$ is a valuation extending $v\in\sum_K\setminus S$, then $\rho(F_w)$ is well defined. Since $\theta$ is algebraic, the eigenvalues of $\theta(\rho(F_w))$ in $\mathrm{GL}_m(\mathbb{Q_\ell})$ are also algebraic numbers. Indeed, if we write $\rho(F_w)=g_{ss}g_u$ by Jordan decomposition, then the eigenvalues of $\rho(F_w)$ are the same as the eigenvalues of $g_{ss}$ and $\theta(\rho(F_w))=\theta(g_{ss})$ because $G_\ell/[G_\ell,G_\ell]$ is a torus. $g_{ss}$ is contained in some maximal torus in $G_\ell$ and we see that the eigenvalues of $\theta(g_{ss})$ are products of integral powers of eigenvalues of $g_{ss}$; therefore the eigenvalues of $\theta(\rho(F_w))$ are algebraic.\\

\noindent\textbf{Proposition 3.18.} Let $\rho:G_K\rightarrow \mathrm{GL}_n(\mathbb{Q}_\ell)$ be a semisimple, rational $\ell$-adic representation of $K$, $G_\ell$ the Zariski closure of $\rho(G_K)$ in $\mathrm{GL}_{n,\mathbb{Q}_\ell}$ and $K^\mathrm{conn}$ the field corresponding to $G_\ell^\circ$. Then $G_\ell$ is reductive, the dimension of the center of $G_\ell^\circ$ is less than or equal to $d_{K^\mathrm{conn}}$. \\

\noindent\textbf{Proof.}  The algebraic group $G_\ell$ is reductive with Lie algebra $\mathfrak{g}_\ell$. We just need to estimate the dimension of the center of $G_\ell^\circ$. We may assume $G_\ell$ is connected by Theorem 1.3. Consider the composition of maps where $\theta$ is defined above,
\begin{equation*}
\theta\circ\rho: G_K\rightarrow G_\ell(\mathbb{Q_\ell})\rightarrow G_\ell/[G_\ell,G_\ell](\mathbb{Q_\ell})\rightarrow \mathrm{GL}_m(\mathbb{Q_\ell}).
\end{equation*}
The quotient map $G_K\rightarrow G_K^\mathrm{ab}$ factors through the composition, hence induces an abelian $\ell$-adic representation 
\begin{equation*}
\Psi: G_K^\mathrm{ab}\longrightarrow \mathrm{GL}_m(\mathbb{Q_\ell}).
\end{equation*}
\noindent The discussion above implies the eigenvalues of $\Psi(F_w)$ are algebraic. Therefore, there exists some integer $N$ such that $(\Psi)^N$ is locally algebraic by Remark 3.15.\\

\noindent Apply Theorem 3.13 to $(\Psi)^N: G_K^\mathrm{ab}\rightarrow \mathrm{GL}_m(\mathbb{Q_\ell})$, we have 
\begin{equation*}
(\Psi)^N=\phi\circ\epsilon_\ell. 
\end{equation*}

\noindent We see that $\phi(S_\mathfrak{m}(\mathbb{Q}_\ell))$ contains the image of $(\Psi)^N$. Since $\Psi(G_K^\mathrm{ab})$ is abelian, the dimensions of $(\Psi)^N(G_K^\mathrm{ab})$ and $\Psi(G_K^\mathrm{ab})$ as $\ell$-adic Lie groups are the same. This implies that $(\Psi)^N(G_K^\mathrm{ab})$ is also Zariski dense in $j(G_\ell/[G_\ell,G_\ell])$. Since $S_\mathfrak{m}(\mathbb Q_\ell)$ is Zariski dense in $S_\mathfrak{m}$ by Theorem 3.8(2), we obtain 
\begin{equation*}
j(G_\ell/[G_\ell,G_\ell]) \subset \phi(S_\mathfrak{m})
\end{equation*}
and conclude that the dimension of the center of $G_\ell$ is less than or equal to the dimension of $S_\mathfrak{m}$ which is $d_K$
by Theorem 3.8(4).\qed\\

\noindent\textbf{Theorem 3.19.} Let $\{\rho_\ell\}_{\ell\in \mathscr{P}}$ be a semisimple, compatible system of $\ell$-adic representations of a number field $K$. 
\begin{equation*}
\rho_\ell:G_K\longrightarrow \mathrm{GL}_n(\mathbb{Q}_\ell).
\end{equation*}
Let $T_\ell$ be a maximal torus of $G_\ell^\circ$ and denote the embedding of $T_\ell$ into $\mathrm{GL}_{n,\mathbb{Q}_\ell}$ by $\Phi_\ell$. Then the triple $(([G_\ell^\circ,G_\ell^\circ]\cap (T_\ell))^\circ,T_\ell,\Phi_\ell)$ is independent of $\ell$. Therefore, the formal character of the tautological representation $(G_\ell^\circ)^\mathrm{der}\hookrightarrow \mathrm{GL}_{n,\mathbb{Q}_\ell}$ and hence the semisimple rank of $G_\ell^\circ$ are independent of $\ell$.\\

\noindent\textbf{Proof.} Assume $K=K^\mathrm{conn}$. Choose some Serre group $S_\mathfrak{m}$ associated to $K$, then it induces an abelian, semisimple, compatible system $\{\epsilon_\ell\}_{\ell\in \mathscr{P}}$ of $\ell$-adic representations with values in $S_\mathfrak{m}$.
\begin{equation*}
\epsilon_\ell:G_K^\mathrm{ab}\longrightarrow S_\mathfrak{m}(\mathbb{Q}_\ell).
\end{equation*}
The image of $\epsilon_\ell$ is an open subgroup of $S_\mathfrak{m}(\mathbb{Q}_\ell)$. Now we choose some faithful representation  $i:S_\mathfrak{m}\rightarrow \mathrm{GL}_{m,\mathbb{Q}}$ over $\mathbb{Q}$, then by base change with $\mathbb{Q}_\ell$ and composing $\epsilon_\ell$ with $i$ we get an abelian, semisimple, compatible system of $\ell$-adic representations of $K$. Still denote it by $\{\epsilon_\ell\}_{\ell\in \mathscr{P}}$:
\begin{equation*}
\epsilon_\ell:G_K^\mathrm{ab}\longrightarrow \mathrm{GL}_m(\mathbb{Q}_\ell).
\end{equation*}
Now consider the system of $\ell$-adic representations $\{\beta_\ell:=\rho_\ell\oplus\epsilon_\ell\}_{\ell\in \mathscr{P}}$, 
\begin{equation*}
\beta_\ell:G_K\longrightarrow \mathrm{GL}_{n+m}(\mathbb{Q}_\ell),
\end{equation*}
which is semisimple and compatible and denote the Zariski closure of the image by $B_\ell$. Since $[\beta_\ell(G_K),\beta_\ell(G_K)]=[\rho_\ell(G_K),\rho_\ell(G_K)]\times \{\mathrm{Id}_m\}$, the Zariski closure of $[\beta_\ell(G_K),\beta_\ell(G_K)]$ is isomorphic to the Zariski closure of $[\rho_\ell(G_K),\rho_\ell(G_K)]$ in $G_\ell$. On the other hand, $\rho_\ell(G_K)\times\rho_\ell(G_K)$ is dense in $G_\ell\times G_\ell$ and the image of the commutator morphism 
\begin{equation*}
[~,~]:G_\ell\times G_\ell\rightarrow G_\ell
\end{equation*}
is closed because $G_\ell$ is connected \cite[Chap. 1 $\mathsection2.3$]{Borel}. We conclude that the Zariski closure of $[\beta_\ell(G_K),\beta_\ell(G_K)]$ is $[G_\ell,G_\ell]\times\{\mathrm{Id}_m\}$. The group $[G_\ell,G_\ell]$ is embedded as a closed normal  subgroup of $B_\ell\subset G_\ell\times S_\mathfrak{m}$. So we have \begin{equation*}
B_\ell/[G_\ell,G_\ell]\subset G_\ell/[G_\ell,G_\ell]\times S_\mathfrak{m}.\\
\end{equation*}

\noindent Choose an embedding $G_\ell/[G_\ell,G_\ell]\hookrightarrow \mathrm{GL}_{n',\mathbb{Q}_\ell}$. Denote the following map by $\alpha_\ell$,
\begin{equation*}
\alpha_\ell:G_K\rightarrow B_\ell(\mathbb Q_\ell)\rightarrow (B_\ell/[G_\ell,G_\ell])(\mathbb Q_\ell)\hookrightarrow \mathrm{GL}_{n'}(\mathbb{Q}_\ell)\times\mathrm{GL}_{m}(\mathbb{Q}_\ell)\subset\mathrm{GL}_{n'+m}(\mathbb{Q}_\ell).
\end{equation*}
 This map is semisimple and factors through $G_K^\mathrm{ab}$. Denote the Zariski closure of $\alpha_\ell(G_K)$ by $C_\ell$ which is diagonalizable. The eigenvalues of the Frobenius elements whenever defined are algebraic. Therefore $\alpha_\ell$ is almost locally algebraic by Remark 3.15, i.e. $\alpha_\ell^N$ is locally algebraic for some positive integer $N$. Since the dimensions of the abelian image of $\alpha_\ell$ and $\alpha_\ell^N$ are equal, we have 
\begin{equation*}
\mathrm{dim}(C_\ell)\leq  d_K
\end{equation*}
by the same argument in the last paragraph of the proof of Proposition 3.18. Observe that the projection of $C_\ell$ to the second factor is $i(S_\mathfrak{m})$, we obtain 
\begin{equation*}
\mathrm{dim}(C_\ell)\geq \mathrm{dim}(S_\mathfrak{m})=d_K.
\end{equation*}
We conclude that $\mathrm{dim}(C_\ell)=  d_K$ is independent of $\ell$ by Theorem 3.8(4). Since $C_\ell$ is isomorphic to $B_\ell/[G_\ell,G_\ell]$, we have an exact sequence of reductive algebraic groups for each $\ell$
\begin{equation*}
0\rightarrow [G_\ell,G_\ell]\rightarrow B_\ell\rightarrow C_\ell\rightarrow 0.
\end{equation*}

\noindent Finally, since $\{\beta_\ell\}_{\ell\in P}$ is a semisimple, compatible system of $\ell$-adic representations of $K$, the rank of $B_\ell$ is independent of $\ell$ by Theorem 1.4. Together with the $\ell$-independence of $\mathrm{dim}(C_\ell)$ and the exact sequence above, we obtain $\ell$-independence of the rank of $[G_\ell,G_\ell]$. Therefore, the dimension of the center of $G_\ell$ is also independent of $\ell$ by Theorem 1.4 again.\\

\noindent Let $\pi_1$ and $\pi_2$ be the projection to the first $n$ coordinates and the last $m$ coordinates respectively.  By base change with $\mathbb{C}$, we assume $T_\ell$ and $\epsilon_\ell(G_K)$ are diagonalized. Also denote the Zariski closure of $\beta_\ell(G_K)$ in $\mathrm{GL}_{m+n,\mathbb{C}}$ by $B_\ell$. We know that the semisimple ranks of $G_\ell$ and $B_\ell^\circ$ are equal by the previous paragraph. Let $D^{n+m}$ be the group of diagonal matrices in $\mathrm{GL}_{n+m,\mathbb{C}}$. Let $D_\ell$ be the connected component of  $B_\ell\cap D^{n+m}$. It is a maximal torus of $B_\ell$. Indeed, $\pi_2(D_\ell)$ has dimension equal to dim($S_\mathfrak{m}$) and $D_\ell$ contains the connected component of $[B_\ell^\circ,B_\ell^\circ]\cap D^{n+m}$ which is equal to the connected component of $([G_\ell,G_\ell]\times\{\mathrm{Id}_m\})\cap D^{n+m}$ having dimension equal to the dimension of $[G_\ell,G_\ell]\cap T_\ell$. Since $\pi_2(([G_\ell,G_\ell]\times\{\mathrm{Id}_m\})\cap D^{n+m})$ is trivial, this implies the dimension of $D_\ell$ is equal to the rank of $B_\ell$. In other words, we could pick for each $\ell$ a diagonalized maximal torus $D_\ell$ of $B_\ell$ such that $\pi_1(D_\ell)=T_\ell$. Since the systems $\{\rho_\ell\}_{\ell\in \mathscr{P}}$ and $\{\epsilon_\ell\}_{\ell\in\mathscr{P}}$ are both compatible, every $D_\ell$ is conjugate by a permutation which permutes the first $n$ coordinates and the last $m$ coordinates. Therefore, we may assume $\Theta_\ell: D_\ell\hookrightarrow \mathrm{GL}_{n+m,\mathbb{C}}$ is independent of $\ell$ and write $\Theta_\ell:=\Theta$. 
\begin{equation*}
\Theta:=(\phi_1,...,\phi_n,\phi_{n+1},...,\phi_{n+m}).
\end{equation*}
Since $[B_\ell^\circ,B_\ell^\circ]\cap D_\ell= ([G_\ell,G_\ell]\times\{\mathrm{Id}_m\})\cap D_\ell   \subset$ Ker$(\pi_2)$ and $\pi_2$ is an isogeny of the center of $B_\ell^\circ$ to its image, we have 
\begin{equation*}
([B_\ell^\circ,B_\ell^\circ]\cap D_\ell)^\circ = \mathrm{Ker}(\pi_2|_{D_\ell})^\circ
\end{equation*}
which is independent of $\ell$. By projecting on the first $n$ coordinates, we have
\begin{equation*}
([G_\ell,G_\ell]\cap T_\ell)^\circ = \pi_1(([B_\ell^\circ,B_\ell^\circ]\cap D_\ell)^\circ) = \pi_1(\mathrm{Ker}(\pi_2|_{D_\ell})^\circ)
\end{equation*}
is also independent of $\ell$. Hence if we denote $(\phi_1,...,\phi_n)$ by $\Phi_\ell$, then the triple \begin{equation*}
(([G_\ell,G_\ell]\cap T_\ell)^\circ, T_\ell,\Phi_\ell)
\end{equation*}
 is independent of $\ell$. It follows that the formal character of $(G_\ell^\circ)^\mathrm{der}\hookrightarrow \mathrm{GL}_{n,\mathbb{Q}_\ell}$ is independent of $\ell$, where $(G_\ell^\circ)^\mathrm{der}:=[G_\ell^\circ,G_\ell^\circ]$.\qed\\

\noindent(\textbf{3.20})  Let's focus on the level of Lie algebra. By Theorem 3.19,  we have a pair $((\mathfrak{g}_\ell)_{ss}\otimes\mathbb{C}, \Phi_\ell)$  for each $\ell$, where $(\mathfrak{g}_\ell)_{ss}\otimes\mathbb{C}$ is a complex semisimple Lie algebra and $\Phi_\ell$ is a faithful representation of $(\mathfrak{g}_\ell)_{ss}\otimes\mathbb{C}$ to an $n$-dimensional complex vector space such that 
\begin{enumerate}
	\item The rank of $(\mathfrak{g}_\ell)_{ss}\otimes\mathbb{C}$ is independent of $\ell$.
	\item The formal character of $\Phi_\ell$ is independent of $\ell$ (see (\textbf{2.1})).
\end{enumerate}
\noindent By Theorem 2.14, 2.17, 2.19, and Proposition 2.20, we obtain the following theorem.\\

\noindent \textbf{Theorem 3.21} Let $\{\rho_\ell\}_{\ell\in \mathscr{P}}$ be a semisimple, compatible system of $\ell$-adic representations of a number field $K$. 
\begin{equation*}
\rho_\ell:G_K\longrightarrow \mathrm{GL}_n(\mathbb{Q}_\ell).
\end{equation*}
Then  the semisimple parts of $\mathfrak{g}_\ell\otimes\mathbb{C}$ satisfy equal-rank subalgebra equivalence (Definition 2.18) for all $\ell$ which is equivalent to the number of $A_n$ factors for $n\in\{6,9,10,...\}$ and the parity of $A_4$ factors of $\mathfrak{g}_\ell\otimes \mathbb{C}$ are independent of $\ell$.\\

\noindent \textbf{Remark 3.22} Let $\{\rho_\lambda\}$ be a compatible system of semisimple, $E_\lambda$-adic representations of $G_K$ \cite[Chap. 1 $\mathsection2.3$]{Serre3}. Locally algebraicity is also defined for abelian $\lambda$-adic representation \cite[$\mathsection2$]{Hen}. Since Theorem 1.3, 1.4, 3.13, 3.14 (see \cite[Thm. 2]{Hen}), and Remark 3.15 still hold analogously for $E_\lambda$-adic representations (the morphisms are then defined over $E_\lambda$), one can prove that Proposition 3.18, Theorem 3.19, and Theorem 3.21 are also true in $E_\lambda$-adic case using identical arguments.\\

\begin{center}
\textbf{$\mathsection4$ Abelian varieties and Galois representations}
\end{center}

Let $A$ be an abelian variety of dimension $g$ over a field $K$, finitely generated over $\mathbb{Q}$. Let $\rho_\ell$ denote the action of $G_K$ on
\begin{equation*}
V_\ell(A):=(\lim_{\leftarrow}A[\ell^n])\otimes\mathbb{Q}_\ell\cong \mathbb{Q}_\ell^{2g}.
\end{equation*}
The image of $\rho_\ell$ is an $\ell$-adic Lie group; denote its Lie algebra by $\mathfrak{g}_\ell$.\\

There exists an abelian scheme $E$ over a smooth variety $X$ defined over a number field $k$ such that the function field of $X$ is $K$ and $E_\eta=A$ where $\eta$ is the generic point of $X$ (see, e.g. Milne \cite[$\mathsection20$]{Mil}). Every closed point $x$ of $X$ induces an $\ell$-adic representation of $\textbf{k}(x)$ given by the Galois action of $G_{\textbf{k}(x)}$ on the $\ell$-adic Tate module of $E_x$, here $\textbf{k}(x)$ is the residue field of $x$ which is a finite extension of $k$.
\begin{equation*}
(\rho_\ell)_x:G_{\textbf{k}(x)}\longrightarrow \mathrm{GL}(V_\ell(E_x))
\end{equation*}
Denote the Lie algebra of the image of $(\rho_\ell)_x$ by $(\mathfrak{g}_\ell)_x$, we have $(\mathfrak{g}_\ell)_x\subset\mathfrak{g}_\ell$ by specialization (see Hui \cite[$\mathsection1$]{Hui}). By \cite[$\mathsection1$]{Serre1}, there always exists a closed point $x\in X$ such that $(\mathfrak{g}_\ell)_x= \mathfrak{g}_\ell$. Therefore, we have $(\mathfrak{g}_\ell)_x= \mathfrak{g}_\ell$ for any prime $\ell$ by \cite[Thm. 2.5]{Hui}. Since the system $\{(\rho_\ell)_x\}_{\ell\in P}$ of $\ell$-adic representation of $\textbf{k}(x)$ is compatible ($\textbf{k}(x)$ is a number field) and semisimple (Faltings \cite{F}), the semisimple parts of $(\mathfrak{g}_\ell)_x\otimes\mathbb{C}=\mathfrak{g}_\ell\otimes\mathbb{C}$ satisfy equal-rank subalgebra equivalence in  Definition 2.18  by Theorem 3.21. Thus, we obtain the following result.\\

\noindent \textbf{Theorem 4.1} Let $A$ be an abelian variety of dimension $g$ over a field $K$, finitely generated over $\mathbb{Q}$. Then we have the following system $\{\rho_\ell\}_{\ell\in \mathscr{P}}$ of $\ell$-adic representation of $K$
\begin{equation*}
\rho_\ell: G_K\longrightarrow \mathrm{GL}(V_\ell(A)).
\end{equation*}
Then the semisimple parts of $\mathfrak{g}_\ell\otimes\mathbb{C}$ satisfy equal-rank subalgebra equivalence (Definition 2.18) for all $\ell$ which is equivalent to the number of $A_n$ factors for $n\in\{6,9,10,...\}$ and the parity of $A_4$ factors of $\mathfrak{g}_\ell\otimes \mathbb{C}$ are independent of $\ell$.\\

\section*{Acknowledgments} This paper contains the main results of my PhD thesis. I am grateful to my advisor, Michael Larsen for many useful conversations during the course of this work and helpful comments on the exposition of the paper. I am also indebted to his generosity, guidance, and continuous support throughout these years. I would also like to thank the Department of Mathematics at Indiana University for providing excellent environment for studies and research.\\

\end{document}